\def\citen#1{\if@filesw \immediate\write \@auxout {\string\citation{#1}}\fi%
\@tempcntb\m@ne \let\@h@ld\relax \def\@citea{}%
\@for \@citeb:=#1\do {\@ifundefined {b@\@citeb}%
    {\@h@ld\@citea\@tempcntb\m@ne{\bf ?}%
    \@warning {Citation `\@citeb ' on page \thepage \space undefined}}%
    {\@tempcnta\@tempcntb \advance\@tempcnta\@ne
    \setbox\z@\hbox\bgroup\ifcat0\csname b@\@citeb \endcsname \relax
    \egroup \@tempcntb\number\csname b@\@citeb \endcsname \relax
    \else \egroup \@tempcntb\m@ne \fi \ifnum\@tempcnta=\@tempcntb
    \ifx\@h@ld\relax \edef \@h@ld{\@citea\csname b@\@citeb\endcsname}%
    \else \edef\@h@ld{\hbox{--}\penalty\@highpenalty
    \csname b@\@citeb\endcsname}\fi
    \else \@h@ld\@citea\csname b@\@citeb \endcsname \let\@h@ld\relax \fi}%
\def\@citea{,\penalty\@highpenalty\hskip.13em plus.13em minus.13em}}\@h@ld}
\def\@citex[#1]#2{\@cite{\citen{#2}}{#1}}%
\def\@cite#1#2{\leavevmode\unskip\ifnum\lastpenalty=\z@\penalty\@highpenalty\fi%
  \ [{\multiply\@highpenalty 3 #1%
  \if@tempswa,\penalty\@highpenalty\ #2\fi}]}   %
\makeatother 

\catcode`\@=12

\def\ee           {\mathbf{e}}

\def\xx           {\mathbf{x}}
\def\yy           {\mathbf{y}}

\def\al           {\alpha}

\def\ga           {\gamma}
\def\de           {\delta}

\def\vep          {\varepsilon}

\def\ka           {\kappa}

\def\om           {\omega}
\def\si           {\sigma}

\def\DE           {\Delta}
\def\LA           {\Lambda}
\def\OM           {\Omega}
\def\PI           {\Pi}
\def\PHI          {\Phi}
\def\PSI          {\Psi}

\def\calc  {{\cal C}}

\def\cale  {{\cal E}}

\def\calh  {{\cal H}}

\def\calo  {{\cal O}}

\def\cals  {{\cal S}}

\def\ssh   {\mathsf{h}}

\def\sshh  {\mathsf{H}}
\def\ssll  {\mathsf{L}}
\def\ssuu  {\mathsf{U}}

\def\sshol {\mathsf{hol}}

\def\id           {\mathbf{1}}

\def\complex       {\mathbb{C}}

\def\reals         {\mathbb{R}}
\def\zet           {\mathbb{Z}}

\def\frag          {\mathfrak{g}}

\def\ra   {\rightarrow}

\def\lra  {\longrightarrow}

\def\lmt  {\longmapsto}

\def\cal           {current algebra}

\def\tft           {topological field theory}

\def\GLN           {GL(n,\complex)}

\def\gln           {\mathfrak{gl}(n,\complex)}

\long\def\query#1{\hskip 0pt{\vadjust{\everypar={}\small\vtop to 0pt{\hbox{}%
     \vskip -13pt\rlap{\hbox to 49.0pc{\hfil{\vtop{\hsize=8pc\tolerance=6000%
     \hfuzz=.5pc\rightskip=0pt plus 3em\noindent#1}}}}\vss}}}}%

\def\futnote#1     {\footnote{~#1}\ }
\long\def\labl#1   {\label{#1}\ee \ifnum\draftcontrol=1
                   \mbox{ }\\[-12 mm]\query{#1}\\[5 mm] \fi}
\long\def\Labl#1#2 {\label{#1#2}\ee\ifnum\draftcontrol=1
                   \mbox{ }\\[-12 mm]\query{#1#2}\\[5 mm] \fi}

\def\ddcl {\frac{\stackrel{\leftarrow}{\delta}}{\delta C}  }
\def\ddcr {\frac{\stackrel{\rightarrow}{\delta}}{\delta C} }
\def\td   {d\!\!^-}
\def\ori  {\mathrm{or}}

\documentclass[12pt]{article}
\usepackage{amssymb,amsfonts,calc,psfig,epsfig,psfrag,latexsym,wasysym}
\usepackage{amsthm,ifthen}
\begin{document}


\begin{center} 

\vskip 15mm

{\Large\bf Topological Field Theory Interpretation\\[5mm]
of String Topology}\\[22mm]
{\large Alberto S. Cattaneo$\;^{2,\star}$ }, \,
{\large J\"urg Fr\"ohlich$\;^{1,\ast}$ }, \,
{\large Bill Pedrini$\;^{1,\dag}$ } , \, \\[5mm]
$^1\;$ Institut f\"ur Theoretische Physik \\
ETH H\"onggerberg\\ CH\,--\,8093\, Z\"urich\\[5mm]
$^2\;$ Institut f\"ur Mathematik \\
Universit\"at Z\"urich\\
Winterthurerstrasse 190\\ CH\,--\,8057\, Z\"urich\\[7mm]
$^{\star}\;$\textsf{asc@math.unizh.ch}\\[2mm]
$^{\ast}\;$\textsf{juerg@itp.phys.ethz.ch}\\[2mm]
$^{\dag}\;$\textsf{pedrini@itp.phys.ethz.ch}
\end{center}

\vskip 25mm

\begin{abstract}
\noindent
The string bracket introduced by Chas and Sullivan is reinterpreted from
the point
of view of topological field theories in the Batalin--Vilkovisky
or BRST formalisms.
Namely, topological action functionals for gauge fields (generalizing
Chern--Simons and $BF$ theories)
are considered together with generalized Wilson loops. The latter
generate
a (Poisson or Gerstenhaber) algebra of functionals with values in the
$S^1$-equivariant
cohomology of the loop space of the manifold on which the theory is
defined.  It is proved that, in the case of
$\GLN$ with standard representation, the (Poisson or BV)
bracket of two generalized Wilson
loops applied to two cycles is the same as the generalized Wilson loop
applied to the string bracket of the cycles.
Generalizations to other groups are briefly described.
\end{abstract}
\newpage


\section{Introduction}

In this paper we study the ``string homology'' defined by Chas and Sullivan \cite{chassull} (see also \cite{cj})
and its algebraic structure from the cohomological point of view of
topological field theory (TFT) \cite{sch,blth}.
String homology provides new topological invariants for general, oriented $d$-dimensional manifolds without boundary. The topological field theory underlying our analysis is a {genera-lization} of three-dimensional Chern-Simons theory, \cite{witt}. It can be defined over an arbitrary differentiable, oriented, $d$-dimensional manifold, $M$, without boundary. Its formulation requires the data of a Lie group $G$\ and a connection, $A$, on a principal $G$-bundle, $P$, over $M$. 

In the main body of this paper we focus our attention on the example where $G=\GLN$, $P$\ is the trivial bundle, $P=M\times G$, and where $A$\ is a flat connection on $P$. But, in the last section of this paper, we sketch the necessary extensions of our arguments to cover more general situations.

We shall study the classical version of our ``topological field theory''; but a few remarks on its quantization are contained in the last section.

Our topological field theory is constructed by making use of the Batalin-Vilkovisky {forma-lism}
or the BRST formalism, depending on whether $d$\ is odd or even; see e.g. \cite{ht}. For the convenience of the reader we recall some key features of these formalisms.

The $BV$\ formalism has been invented as a tool to quantize field theories 
in the Lagrangian formalism
with a large (infinite) number of (infinitesimal) symmetries, 
for example gauge theories.
The space, $\calc_0$, of classical field configurations of such a theory is first augmented by introducing \textit{ghosts}, and second by introducing \textit{antifields for fields and ghosts} in equal number as the fields and the ghosts. The extended configuration space, $\calc$, thus obtained can be viewed as an
(odd-symplectic)
supermanifold, the fields, ghosts and antifields for fields and ghosts being local even or odd (Darboux)
coordinates on it. 
The superfunctions on $\calc$\ form the supercommutative 
algebra of ``preobservables'',
denoted by $\calo$. 
This algebra is equipped with a natural $\zet_2$-grading, $|\cdot|$, 
and
is furnished by construction with a non-degenerate, odd bracket, $\{\cdot;\cdot\}$,
\begin{equation}
  \begin{array}{lccc}
  \{\cdot;\cdot\}:&
  \calo\times\calo &\lra& \calo\\
  & (\;O_1 , O_2\;) & \mapsto & \{O_1;O_2\}
  \end{array}
\end{equation}
satisfying graded versions of \textit{antisymmetry}, of the \textit{Leibnitz rule}, and of the \textit{Jacobi identity}. This is equivalent to saying that $(\calo,\{\cdot;\cdot\})$\ is a \textit{Gerstenhaber algebra}.
Choosing local ``Darboux coordinates'', $\phi^a,\phi^{\dag}_a$, on $\calc$, for example interpreting the $\phi^a$'s as ``fields'' (fields and ghosts) and the $\phi^{\dag}_a$'s as ``antifields'' (antifields for fields and 
ghosts),\footnote{$\phi^a$ and $\phi^{\dag}_a$ are assigned opposite Grassmann
parity.} 
the bracket can be expressed as 
\begin{equation}
  \{O_1;O_2\}=
     O_1\frac{\stackrel{\leftarrow}{\partial}}{\partial \phi^a}
      \frac{\stackrel{\rightarrow}{\partial}}{\partial \phi^{\dag}_a})_2
     -
     O_1\frac{\stackrel{\leftarrow}{\partial}}{\partial \phi^{\dag}_a}
      \frac{\stackrel{\rightarrow}{\partial}}{\partial \phi^a}O_2
  \quad.
\end{equation}
In classical theory, one attempts to construct an action functional $S$\ of degree zero satisfying the \textit{classical master equation}
\begin{equation}
  \{S;S\}=0 \quad.
\end{equation}
Such an action functional equips $\calo$\ with the structure of a \textit{differential algebra}. The differential, $\de$, is given by
\begin{equation}\label{eq:diffe}
  \de O=\{S;O\} \quad (O\in\calo)\quad.
\end{equation}
Because the bracket is odd and $|S|=0$,
\begin{equation}
  |\de O|=|O|+1\quad.
\end{equation}
The classical master equation for $S$\ and the graded Jacobi identity imply that $\de$\ is nilpotent, i.e.,
\begin{equation}
  \de^2=0 \quad.
\end{equation}
The cohomology of $\de$, $H^{\ast}_{\de}$, is called the algebra of
``\textit{observables}'' of the theory. Thanks to the graded Leibnitz rule it
is indeed an algebra. The master equation and the graded Jacobi identity can
be used to show that the bracket descends to cohomology, and $H^{\ast}_{\de}$\
thus has the structure of a \textit{Gerstenhaber algebra}.

The structure described above is well suited to formulate a topological field theory yielding the cohomological version of the results of Chas and Sullivan, provided the dimension $d$\ of the underlying manifold $M$\ is \textit{odd}.
When $d$\ is \textit{even} we must actually follow the (Hamiltonian) BRST formalism.
The latter was developed to quantize theories with (first-class)
constraints. The classical phase space, $\calc_0$, is augmented by introducing
\textit{ghosts} and \textit{antighosts} in equal number. The extended space,
$\calc$, thus obtained can be considered as a supermanifold, the fields,
ghosts and antighosts being (even or odd) coordinates on it. 
The algebra, $\calo$, of preobservables is defined to be the algebra of
superfunctions on $\calc$. By construction, $\calo$ is furnished with a
non-degenerate, even bracket. Thus the algebra $\calo$ has the structure of a \textit{super-Poisson algebra}.
The action $S$, now more appropriately called \textit{BRST generator}, is odd
($|S|=1$). The differential $\de$ on the algebra of preobservables is still
defined by (\ref{eq:diffe}), it has degree 1 and is nilpotent. The cohomology $H^{\ast}_{\de}$\ of $\de$\ now has the structure of a \textit{super-Poisson algebra}. (Observe that $H^{0}_{\de}$ describes the algebra of functions on the
reduced phase space, but in general other cohomology groups may be
nontrivial, too.) 

The Lagrangian BV formalism and Hamiltonian BRST (or BFV) 
formalism are related to each other:
after gauge fixing of the BV master action, which requires the elimination 
of the antifields by expressing them as appropriate functions of the fields, 
one finds an action
for which the Legendre transformation to pass to the Hamiltonian formalism can
be pursued; the Hamiltonian so obtained has BRST symmetry, and the BRST
generator can be constructed. 
For more details we refer the reader to Appendix \ref{app:bvbrst}, where the connection between the two formalisms is illustrated for our topological field theory.

In this paper we start directly from an \textit{extended field space} $\calc$\
and a master action (BRST generator) $S$\ satisfying the classical master
equation, see Section \ref{sec:gengf}, without asking whether the theory comes from a classical Lagrangian (or Hamiltonian) theory.

Field configurations of our theory are differential forms, $C$, on $M$ with values in the tensor product of a supercommutative algebra, $\cale$, 
with the metric
\footnote{A Lie algebra endowed with a 
non-degenerate, Ad-invariant inner product is called metric. In particular,
semi-simple Lie algebras with the Killing form are metric. But so are
abelian Lie algebras with any non-degenerate inner product.}
Lie algebra $\frag$\ of the Lie group $G$. 
For simplicity, we suppose that the metric on $\frag$ is given
by the trace in a representation $\rho_0$. 
The forms $C$\ have total degree $|C|=1$, where the mod 2 grading $|\cdot|$\ takes account of both the form degree and the $\cale$-degree.
The space of field configurations, $\calc$, can be considered as a supermanifold with a natural odd (even) bracket; this gives the space of ($\cale$-valued) superfunctions, $\calo$, the structure of a Gerstenhaber (super-Poisson) algebra.
The action functional, $S$, is chosen to be the 
\textit{``Chern-Simons'' action}
\begin{equation}\label{eq:action}
  S[C]=\int_M\mathrm{tr}_{\rho_0}\left[
       \frac{1}{2}Cd_AC+\frac{1}{3}C^3\right]
  \quad,
\end{equation}
where $d_A$\ is the covariant exterior derivative
(w.r.t.\ the flat connection $A$) over $M$. 
Of course, in the integrand of (\ref{eq:action}) only the part of total 
form degree $d$ contributes. 
It is not hard to show that the action is even (odd), $|S|=0$, ($|S|=1$), 
and that it satisfies the master equation, $\{S;S\}=0$.

Observables of these theories can be constructed as follows.
Let $\ssll M$\ denote the space of marked, parametrized loops in $M$. It
carries an obvious circle action. String space, $\cals M$, is defined as the
quotient of $\ssll M$\ by this circle action; see Section
\ref{sec:strsp}. From the connection $A$\ and the forms $C$\ one can
construct, using Chen's iterated integrals (``Dyson series''),
\textit{generalized holonomies}, $\sshol_{A}(C)$, in a fairly obvious way
explained in Section \ref{sec:genhol}. The trace,
$\ssh_{\rho;A}(C)=\mathrm{tr}_{\rho}\sshol_{A}(C)$, also called
\textit{generalized Wilson loop}, then defines a
(generalized) preobservable with values in $\cale\otimes\OM^{\ast}(\cals M)$,
i.e., a differential form on $\cals M$ whose components take values in a
supercommutative algebra $\cale$. If $a$\ represents a cycle in string homology, $\calh_{\ast}M$, as described in \cite{chassull}, then one can pair $a$\ with $\ssh_{\rho;A}(C)$ by integration,
\begin{equation}
  \int_a\ssh_{\rho;A}(C)
  \quad.
\end{equation}
We shall see in Section \ref{sec:genhol} that $\int_a\ssh_{\rho;A}(C)$\ is an \textit{observable} of the theory, i.e., $\de\int_a\ssh_{\rho;A}(C)=0$, for arbitrary $[a]\in\calh_{\ast}M$.

The main result of this paper, proven in Section \ref{sec:homo}, is the following theorem.\\
\textbf{Theorem}. Let $G=\GLN$, $n=1,2,3,\ldots$, and let $\rho$ denote its standard representation (as matrices on $\complex^n$). Let $A$\ be a flat connection on $M\times G$. Then
\begin{equation}
  \left\{\int_a\ssh;\int_{\bar{a}}\ssh\right\}=\int_{\{a;\bar{a}\}}\ssh
  \quad,
\end{equation}
where $\{a;\bar{a}\}$\ is the Chas-Sullivan bracket, see \cite{chassull}, defined on string homology, and $\ssh$\ is a shorthand notation for $\ssh_{\rho;A}(C)$.$\quad\Box$

The definition of the Chas-Sullivan bracket on string homology and some of its properties are explained in Section \ref{sec:bracket}. The special role played by the groups $\GLN$\ is explained in Section \ref{sec:gln}. As sketched in Section \ref{sec:outlook}, more general Lie groups can be accommodated by replacing the string space by a ``space of chord diagrams'' on the manifold $M$.
Section \ref{sec:outlook} also contains a sketch of various other generalizations (e.g. to nontrivial principal $G$-bundles).\\
\textbf{Acknowledgments}. 

B. P. thanks Carletto Rossi for useful discussions
about generalized holonomies.

A. S. C. acknowledges a three-month invitation at Harvard University during
the Fall Term 2001, and thanks Raoul Bott and David Kazhdan for stimulating
discussions.

A. S. C. thanks partial support by SNF Grant No.\ 20-63821.00 .

\section{A TFT with generalized gauge fields}
  \label{sec:gengf}

In this section, we introduce the topological field theories described
in the Introduction in a mathematically precise fashion. We first describe
the space of field configurations, then we introduce algebras of preobservables
and define the bracket between two preobservables, and, finally, we define
an ``action functional`` satisfying the classical master equation.

\subsection{Field configurations}
The field theory is defined over a differentiable,
 oriented, $d$-dimensional manifold
$M$.

Let $P=M\times G$\ be a (for simplicity trivial) principal bundle over $M$\
with structure group $G$. Denote by $\frag$\ the Lie algebra of $G$,
by $\ssuu\frag$\
the corresponding universal enveloping algebra, and by $\ka(\cdot,\cdot)$\ an
invariant bilinear form on $\frag$, which, for notational simplicity, we suppose to be given
by the trace in some representation $\rho_0$:
$\ka(\cdot,\cdot)=tr_{\rho_0}[\cdot\;\cdot]$.

Let $A$\ be a flat connection on $P$, i.e., $A\in\OM^1(M,\frag)$ with
$dA+\frac{1}{2}[A,A]=0$.

We require the following mathematical objects and concepts. A
superalgebra $X$\ (over $\reals$)
is an algebra furnished with a mod 2 grading $|\cdot|$, such 
that, as a vector space, it has the structure
 $X=X_0\oplus X_1$, with $|x_i|=i$\ for $x_i\in
X_i$, and such that $|x_1x_2|=|x_1|+|x_2|$. A superalgebra is supercommutative 
if $x_1x_2=x_2x_1(-1)^{|x_1||x_2|}$.

Next, let $\cale$\ be a supercommutative algebra (e.g. the algebra of
supernumbers \cite{dewitt}). A superalgebra $X$\ is an $\cale$-bimodule if
$\cale$\ acts on $X$ from the left and the right, with 
$\vep x=x\vep(-1)^{|x||\vep|}$\ and $|\vep x|=|\vep|+|x|$, for arbitrary 
$\vep\in\cale$\ and $x\in X$. 
$\cale$\ is clearly an
$\cale$-bimodule.

Any superalgebra $X$\ can be turned into an $\cale$-bimodule
by considering $X_{\cale}=\cale\otimes_{\reals}X$\ and defining
the grading $|\vep\otimes x|=|\vep|+|x|$,
the left action $\vep_1(\vep_2\otimes x)=(\vep_1\vep_2)\otimes x$,
the right action 
$(\vep_2\otimes x)\vep_1=(\vep_1\vep_2)\otimes x(-1)^{|x||\vep_2|}$,
and the product  
$(\vep_1\otimes x_1)(\vep_2\otimes x_2)=\vep_1\vep_2\otimes
x_1x_2(-1)^{|x_1||\vep_2|}$. 
For notational simplicity, one writes $\vep\equiv\vep\otimes\id$,   
$x\equiv\id\otimes x$\ and $\vep x\equiv\vep\otimes x$. 

Given two superalgebras $X_1$\ and $X_2$\ which are $\cale$-bimodules, 
one may define a tensor product bimodule
$X_1\cdot X_2=X_1\otimes_{\cale}X_2$, which becomes a superalgebra by defining
the grading as $|x_1\otimes x_2|=|x_1|+|x_2|$\ and the product as $(x_1\otimes 
x_2)(y_1\otimes y_2)=x_1y_1\otimes x_2y_2(-1)^{|x_2||y_1|}$. For notational
simplicity one writes $x_1\equiv x_1\otimes\id$, $x_2\equiv \id\otimes x_2$\
and $x_1x_2\equiv x_1\otimes x_2$. Clearly one has that $\cale\cdot X=X$.

Let $\calc^G=\OM^{\ast}(M)_{\cale}\cdot\frag_{\cale}$. The space of field
configurations is defined as
\begin{equation}
  \calc^G_1=\{C\in\calc^G\big{|}|C|=1\}
  \quad.
\end{equation}
We note that the components,
$C^a_{\mu_1\ldots\mu_k}(x)\in\cale$, of a field configuration $C\in\calc^G_1$,
 are bosonic for odd $k$\ and
fermionic for even $k$; ($a$\ labels a basis in $\frag$).

\vspace{4mm}
\subsection{Preobservables}
A generalized preobservable is a functional on the space of field
configurations with values in a superalgebra $X$\ which is also an
$\cale$-bimodule; i.e., it is an element of
\begin{equation}
  \calo^G(X)\equiv\OM^{0}(\calc^G_1,X)
\end{equation}
$\calo^G(X)$
is clearly an $\cale$-bimodule, the grading being given by the grading on $X$.
We shall not indicate the group $G$\ if not necessary.
The space of (ordinary) preobservables is $\calo\equiv\calo(\cale)$.
Though not strictly necessary, the concept of generalized preobservables turns 
out to be very convenient in the following.

The (tensor) product 
of two preobservables is defined as a map from
$\calo(X_1)\times\calo(X_2)$\ to $\calo(X_1\cdot X_2)$\
in the obvious way.

\subsection{Bracket between preobservables}
We begin by defining the two operators
\begin{equation}
  \ddcl,\ddcr:
 \calo(X)\lra\calo(X\cdot\OM^{\ast}(M)_{\cale}\cdot\frag_{\cale})
\end{equation}
as follows:

\begin{equation}\label{eq:defder}
  \left.\frac{d}{dt}\right|_{t=0}O(C+t\eta)=
  \int_M\mathrm{tr}_{\rho_0}\left[\eta\ddcr O\right]=
  (-1)^{d(d+|O|)}\int_M\mathrm{tr}_{\rho_0}\left[O\ddcl\eta\right]
  \quad,
\end{equation}
for $O\in\calo(X)$\ and arbitrary $\eta\in\calc_1$. 
The signs are chosen in such a way that 
these two operators act from the left/right as operators of degree $d+1$,
i.e., such that the Leibnitz rules
\begin{eqnarray}\label{eq:ddcleib}
  \ddcr(O_1O_2)&=&(\ddcr O_1)O_2+(-1)^{|O_1|(d+1)}O_1(\ddcr O_2)\quad, \\
  (O_1O_2)\ddcl&=&(-1)^{|O_2|(d+1)}(O_1\ddcl)O_2+O_1(O_2\ddcl)
\end{eqnarray}
hold.
Moreover, one has
\begin{equation}\label{eq:ddcrddcl}
  \ddcr O=(-1)^{(d+1)|O|+1}O\ddcl \quad.
\end{equation}

Next, we define the bracket, $\{\cdot;\cdot\}$, by
\begin{equation}\label{eq:fieldbr}
  \begin{array}{lccc}
  \{\cdot;\cdot\}:&
  \calo(X_1)\times
  \calo(X_2)&\lra&
  \calo(X_1\cdot X_2)\\
  & (\;O_1 , O_2\;) & \mapsto & 
  \{O_1;O_2\}=(-1)^{|O_1|d}
    \int_M \mathrm{tr}_{\rho_0}\left[O_1\ddcl\ddcr O_2\right]
  \end{array}
  \quad.
\end{equation}
The signs are chosen in such a way that, for $d$\ even, $\{\cdot;\cdot\}$\
 is an even bracket, while for
$d$\ odd it is an odd bracket. In fact, $\{\cdot;\cdot\}$\ has the following 
properties:
\begin{itemize}
\item[(1)] \textit{Antisymmetry},
\begin{equation}
  \{O_1;O_2\}=
    -(-1)^{(|O_1|+d)(|O_2|+d)}\{)_2;O_1\}
    \quad,
\end{equation}
a consequence of (\ref{eq:ddcrddcl});
\item[(2)] \textit{Leibnitz rule} 
\begin{equation}
  \label{eq:leibobs}
  \{O_1;O_2O_3\}=\{O_1;O_2\}O_3+(-1)^{|O_2|(|O_1|+d)}O_2\{O_1;O_3\}
  \quad,
\end{equation}
a consequence of (\ref{eq:ddcleib});
\item[(3)] \textit{Jacobi identity} 
\begin{equation}
  \label{eq:jacobs}
  \{O_1;\{O_2;O_3\}\}=\{\{O_1;O_2\};O_3\}+
                      (-1)^{(|O_1|+d)(|O_2|+d)}\{O_2;\{O_1;O_3\}\}
  \quad,
\end{equation}
which can be checked by using (\ref{eq:ddcrddcl}), (\ref{eq:ddcleib}) 
and the definition (\ref{eq:fieldbr}) .

\end{itemize}

We observe that, for a manifold $A$, for multivector fields 
$v_i\in\OM_{\ast}(A)$\ and for generalized
preobservables $O_i\in\calo(\OM^{\ast}(A)_{\cale})$\ the contraction ($\equiv$\
infinitesimal integration of chains with given orientation) 
can be understood as an
operator, $\iota_{v}$\, acting from the left and of degree $|v|$, namely
\begin{equation}\label{eq:contraction1}
  \iota_{v_1}\{O_1;O_2\}=\{\iota_{v_1}O_1;O_2\} \quad,\quad
  \iota_{v_2}\{O_1;O_2\}=(-1)^{|v_2|(d+|O_1|)} \{O_1;\iota_{v_2}O_2\} \quad.
\end{equation}

An explicit calculation on $\calo$\ reveals that
\begin{equation}\label{eq:brcc}
  \{C^a_{\mu_1\ldots\mu_k}(x);C^b_{\mu_{k+1}\ldots\mu_d}(y)\}=
  (-1)^k \de^{(d)}(x-y)\ka^{ab}\vep_{\mu_1\ldots\mu_k\mu_{k+1}\ldots\mu_d}
  \quad.
\end{equation}

\subsection{BRST/BV generator and observables}
We define an ``action'' functional, $S$, by
\begin{equation}
  S[C]=\int_M 
  \mathrm{tr}_{\rho_0}\left[\frac{1}{2} Cd_AC+
    \frac{1}{3}C^3 \right] \in \calo \quad.
\end{equation}
This functional has total degree $d+1$ and 
is constructed so as to satisfy the 
BV/BRST master equation,
\begin{equation}
  \{S;S\}=0 \quad.
\end{equation}
It is thus to be thought of as a \textit{classical master action} in the 
Lagrangian formalism,
for $d$\ odd,
or as a \textit{classical BRST generator}
 in the Hamiltonian formalism, for $d$\ even. 
Being independent of the choice of a metric on $M$, the field theoretical 
model is
called topological
\footnote{There is a sigma-model construction of $S$\ and $\{\cdot;\cdot\}$,
  obtained by considering the fields $C$\ as maps $\PI TM\lra\PI\frag$\ 
(see \cite{aksz}), where $\PI$\ reverses the parity of the fiber 
in a vector bundle}. 
One can check that,
in a situation where $M^{[d+1]}=M^{[d]}\times\reals$, $d$ even,
 $S^{[d]}$\ is the BRST generator corresponding
to
$S^{[d+1]}$\, after gauge fixing; (see Appendix \ref{app:bvbrst}).

$S$\ defines an odd  differential, $\de$, on the algebra of preobservables
by 
\begin{equation}
  \begin{array}{cccc}
  \de:&
  \calo(X) &\lra& \calo(X)\\
  & O & \mapsto & \{S;O\} 
  \end{array}
  \quad.
\end{equation}
We wish to mention another important property of $S$: 
The bracket between $S$ and a field component $C$ is given by
\begin{equation}\label{eq:deltafields}
  \{S;C\}=(-1)^d(d_AC+C^2)
  \quad,
\end{equation}
or, more explicitly,
\begin{equation}
  \{S;C^a_{\mu_1\ldots\mu_k}(x)\}=(-1)^{d+k}(d_AC+C^2)^a_{\mu_1\ldots\mu_k}(x)
\end{equation}
This is a key equation for proving the fundamental identity
(\ref{eq:ddeltah}), below.

The cohomology of $\de$, $H_{\de}^{\ast}$, defines the algebra of
generalized observables of the topological field theory. Because of
(\ref{eq:leibobs}) and (\ref{eq:jacobs}), respectively, product and bracket
descend to {cohomo-logy}; 
the generalized observables thus have the structure of\\
$\Diamond$\ a super-Poisson algebra (even bracket), for $d$\ even,\\
$\Diamond$\ a Gerstenhaber algebra (odd bracket), for $d$\ odd.\\

\section{The String Space of a manifold}
  \label{sec:strsp}

In this section we define the loop space of a manifold, and, subsequently, the string space
as the quotient of the former by a circle action. Moreover, we describe how to 
define local coordinates on loop- and string space.

\vspace{4mm}
One may define the loop space of a manifold $M$\ as
\begin{equation}
  \ssll M=\{\ga(\cdot):S^1\lra M,\;\
  \ga\;\mathrm{piecewise}\;\mathrm{differentiable}\} 
  \quad.
\end{equation}
Observe that $S^1$\ has a marked point, $0$, if we interpret
$S^1$\ as $\reals/\zet$. 
Therefore a loop can be thought of as a parametrized closed curve in $M$\
with a marked point and a tangent vector in almost every point, 
the parameter $t$ ranging from $0$\ to $1$.

Let $(x^{\mu})_{\mu=1\ldots d}$\ be local coordinates 
on a coordinate patch $U\subset M$.
Then $(\ga^{\mu}(t))_{\mu=1\ldots d,t\in S^1}$\ are corresponding local
coordinates on the patch $\ssll U\subset\ssll M$. (For loops which extend over
different patches, there is a similar construction of local coordinates;
but it is not needed for the purposes of this paper).

Loop space carries an obvious circle action
\begin{equation}
  \begin{array}{ccc}
    S^1\times\ssll M & \lra & \ssll M \\
    \left(s, \ga(\cdot)\right) & \lmt & \ga(\cdot+s)  
  \end{array}
  \quad.
\end{equation}
The string space, $\cals M$\, is defined as the quotient of $\ssll M$\ by this action
\footnote{The string space is a singular manifold, with singularities arising
 at the constant loops/strings, which are 
 fixed points of the circle action.} 
\begin{equation}\label{eq:circac}
  \begin{array}{ccc}
  S^1 & \hookrightarrow & \ssll M\\
  & & \downarrow \pi_{S^1} \\
  & & \cals M
  \end{array}
  \quad.
\end{equation}
A string can thus be thought of as a closed curve in $M$
with a tangent vector in almost every point.

\begin{figure}[h]
\begin{center}
  \input{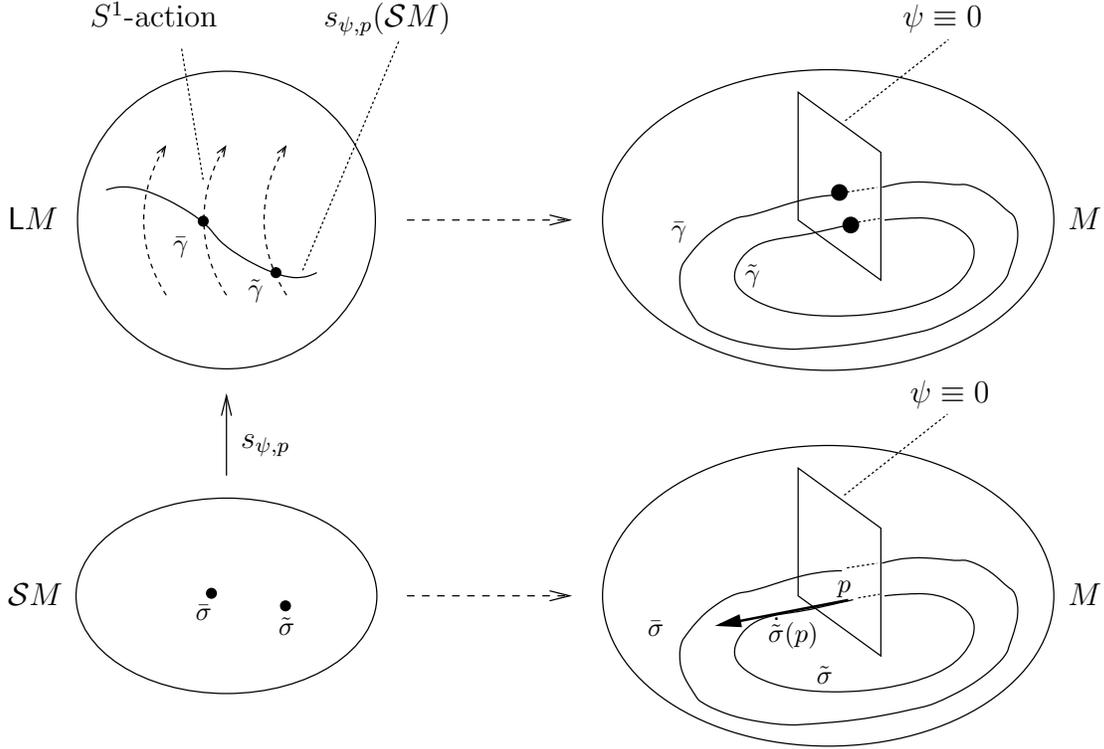}
  \parbox{14cm}{
    \caption{\label{fig:loccoo}Constructing local coordinates on $\cals M$.}}
\end{center}
\end{figure}

Local coordinates on $\cals M$\ can be constructed by choosing a local section
$\cals M\lra\ssll M$\ and then using local coordinates on $\ssll M$; 
see Figure \ref{fig:loccoo}.
More precisely, let $\tilde{\si}\in\cals U\subset\cals M$\ be a nonconstant 
string and
$p$\ a 
point on it such that $\dot{\tilde{\si}}(p)\not=0$.
Let $\psi$\ be a function on $M$\ 
defined in a neighborhood of $p$\ 
such that $\psi(p)=0$\ and
$\big{<}\dot{\tilde{\si}}(p);d\psi(p)\big{>}\neq0$.
A local section $s_{\psi,p}:\cals M\lra\ssll M$\ in a neighborhood of 
$\tilde{\si}$ is uniquely defined by the requirement 
that $\psi(s_{\psi,p}(\bar{\si})(t=0))=0$,
for any string $\bar{\si}$\ that is a 
sufficiently small deformation of $\tilde{\si}$. 
The functions
$(\si^{\mu}(t))_{\mu=1\ldots d,t\in S^1}$, defined as
$\si^{\mu}(t)=\ga^{\mu}(t)\circ s_{\psi,p}$, are then local coordinates on
$\cals M$\ in a neighborhood of $\tilde{\si}$.

We denote by $\calh_{\ast}M$\ the
string homology, properly defined as the $S^1$-equivariant loop space
homology. 
We denote by $\td$ the differential on both loop- and string space.

\section{Generalized holonomies and Wilson loops}
  \label{sec:genhol}

In this section we define generalized Wilson loops as generalized
observables with values in string cohomology.
As such, they can be paired with cycles in string homology,
yielding observables of the \tft.

\vspace{4mm}
\noindent 
We introduce standard simplices
$\left.\DE_n\right|_{t_i}^{t_f}=\{(t_1,\ldots,t_n)\in\reals^n|t_i\leq
t_1\leq\ldots\leq t_n\leq t_f\}$, $\DE_n=\left.\DE_n\right|_{0}^{1}$,
and define the evaluation maps
\begin{equation}
  \begin{array}{ccccc}
  \mathrm{ev}_{n,k}: & \DE_n\times\ssll M & \lra & M & \\
  & \left(t_1,\ldots,t_n;\ga\right) & \lmt & \ga(t_k) &\qquad 1\leq k\leq n
  \end{array}
\end{equation}
The $n$-th order generalized parallel transporter is given by
\begin{equation}
  \left.\sshol^n_{A}(C)\right|_{t_i}^{t_f}=
      \int_{\left.\DE_n\right|_{t_i}^{t_f}}
      \left(
      \left.\mathrm{hol}_{A}\right|_{t_i}^{t_1}\;
      \mathrm{ev}^{\ast}_{n,1}C\;
      \left.\mathrm{hol}_{A}\right|_{t_1}^{t_2}\ldots
      \left.\mathrm{hol}_{A}\right|_{t_{n-1}}^{t_n}\;
      \mathrm{ev}^{\ast}_{n,n}C\;
      \left.\mathrm{hol}_{A}\right|_{t_n}^{t_f}
      \right)
   \quad.
\end{equation}
In this definition the parallel transporter,
$\left.\mathrm{hol}_{A}\right|_{t_k}^{t_{k+1}}=
P\exp\int_{t_k}^{t_{k+1}}\iota_{\dot{\ga}(t)}A$,
of the flat connection $A$ is a function $\DE_n\times\ssll
M\lra\ssuu\frag_{\cale}$; ($P$ denotes path ordering). 
For an expression in local coordinates, see Appendix \ref{app:lochol}.

Thus,$\left.\sshol^n_{A}\right|_{t_i}^{t_f}$ is an element of
$\calo(\OM^{\ast}(\ssll M)_{\cale}\cdot\ssuu\frag_{\cale})$. 
We define generalized parallel transporters,
$\left.\sshol_{A}\right|_{t_i}^{t_f}$, by
\begin{equation}\label{eq:parthol}
  \left.\sshol_{A}(C)\right|_{t_i}^{t_f}=
  \sum_{n=0}^{\infty}\left.\sshol^n_{A}(C)\right|_{t_i}^{t_f}\quad,
\end{equation}
and generalized holonomies by
\begin{equation}
  \sshol_{A}(C)=
  \left.\sshol_{A}(C)\right|_{0}^{1}
  \quad.
\end{equation}
Furthermore, generalized ``Wilson loops'' in a representation $\rho$\ are
defined by
\begin{equation}
  \ssh_{\rho;A}(C)=\mathrm{tr}_{\rho}\sshol_{A}(C)\quad.
\end{equation}
It is worth remarking that the degree of generalized parallel transporters and
generalized Wilson loops is zero, i.e.,
\begin{equation} \label{eq:deghol}
  \left|\sshol_{A}\right|=\left|\ssh_{\rho;A}\right|=0
  \quad.
\end{equation}
Under a gauge transformation, $g:M\lra G$, one finds that
\begin{equation}
  \sshol_{A}(C)=g^{-1}\sshol_{g(A+d)g^{-1}}(gCg^{-1})g\quad,\quad
  \ssh_{\rho;A}(C)=\ssh_{\rho;g(A+d)g^{-1}}(gCg^{-1}) \quad.
\end{equation}

The tangent vectors, $\dot{\ga}$, that generate the circle action on $\ssll M$
define a section of $T\ssll M$. The contraction
  $\iota_{\dot{\ga}}\ssh_{\rho;A}$\ clearly vanishes. Moreover, one finds
  \cite{accr1} that
\begin{equation}\label{eq:dhol}
  \td\ssh_{\rho;A}=\int_0^1d\tau \; \mathrm{tr}_{\rho}\left[
        \left.\sshol_{A}(C)\right|_{0}^{\tau}
        \iota_{\dot{\ga}}
        \mathrm{ev}^{\ast}_{\tau}(d_AC+C^2)
        \left.\sshol_{A}(C)\right|_{\tau}^{1}\right]
  \quad,
\end{equation}
where $\mathrm{ev}_{\tau}:\ssll M\lra M,\;\ga\lmt\ga(\tau)$. 
This implies that the Lie derivative 
$L_{\dot{\ga}}\ssh_{\rho;A}=\iota_{\dot{\ga}}\td\ssh_{\rho;A}$\ vanishes, too. The form 
$\ssh_{\rho;A}$\ is thus horizontal and invariant with respect to 
the circle action,
and thus defines a form on string space.

Comparing (\ref{eq:dhol}) and (\ref{eq:deltafields}), 
we find the fundamental identity \cite{accr1}\cite{accr2}
\begin{equation}\label{eq:ddeltah}
  ((-1)^d\de+\td)\ssh_{\rho;A}=0
  \quad,
\end{equation}
which implies that the trace of the generalized holonomy is an
observable with values in string cohomology
\footnote{There are no problems connected with the singularities 
 of string space,
  since the form vanishes at constant strings.},
\begin{equation}
  \ssh_{\rho;A}\in
  H^{\ast}_{\de}\calo(\calh^{\ast}M)
  \quad,
\end{equation}
and, for a cycle $a\in\calh_{\ast}M$\ in string homology, the pairing
\begin{equation}
  \big{<}a,\ssh_{\rho;A}\big{>}:=\int_a\ssh_{\rho;A}\;\in
  H^{\ast}_{\de}\calo
\end{equation}
defines an observable.

\section{The String Bracket}
  \label{sec:bracket}

In this section we recall how to define a bracket
\begin{equation}
  \{\cdot;\cdot\}:
  \calh_{\ast}M\times\calh_{\ast}M\lra\calh_{\ast}M
\end{equation}
on string homology. 
This definition is taken from the article of Chas and Sullivan
\cite{chassull}, but we give a slightly simplified exposition.

\vspace{4mm}
Define $\cals M^{\times}\subset\cals M\times\cals M$\ as the
space of pairs of strings which intersect transversally at at least one point.
This space is a cycle of codimension $d-2$, with $n-1$-fold self intersections
when the two strings intersect $n$\ times. 
We propose to construct the current corresponding to $\cals M^{\times}$.
The $d$-form
\begin{equation}
  \om^{\times}=\de(x^1-\bar{x}^1)\ldots\de(x^d-\bar{x}^d)
               (dx^1-d\bar{x}^1)\ldots(dx^d-d\bar{x}^d)
  \in\OM^{d}(M\times M)
\end{equation}
is the current for the diagonal in $M\times M$. We define 
\begin{equation}
  C^{\times}=\int_{S^1\times\bar{S}^1}
             (\mathrm{ev}^{\ast}_{1,1}\times\bar{\mathrm{ev}}^{\ast}_{1,1})
             \om^{\times}
  \quad,
\end{equation}
which is a $(d-2)$-current on $\ssll M\times \ssll M$. It is closed, since $\om^{\times}$\ is closed, and the integration domain, $S^1\times\bar{S}^1$, in the above formula has no boundaries.
In local coordinates, it reads
\begin{eqnarray} 
  C^{\times} & = &
    \sum_{k=1}^{d-1}\frac{(-1)^{d+1}}{(k-1)!(d-k-1)!}
    \int_{s=0}^{s=1}ds\int_{\bar{s}=0}^{\bar{s}=1}d\bar{s}
    \; \de^{(d)}(\ga(s)-\bar{\ga}(\bar{s}))
    \; \vep_{\nu_1\nu_2\ldots\nu_k
      \bar{\nu}_{k+1}\bar{\nu}_{k+2}\ldots\bar{\nu}_d}
    \nonumber\\
    & & \label{eq:loccurr}
    \dot{\ga}^{\nu_1}(s)\td\ga^{\nu_2}(s)\ldots\td\ga^{\nu_k}(s)
    \dot{\bar{\ga}}^{\bar{\nu}_{k+1}}(\bar{s})
    \td\bar{\ga}^{\bar{\nu}_{k+2}}(\bar{s})
    \ldots\td\bar{\ga}^{\bar{\nu}_d}(\bar{s}) 
  \quad.
\end{eqnarray}
{}From this expression it is easy to see that it is horizontal, and thus also
invariant with respect to the two circle actions on the two factors of $\ssll
M\times \ssll M$. Hence, $C^{\times}$\ defines a closed $(d-2)$-current 
on $\cals M\times\cals M$. 
Let $(\si,\bar{\si})$\ be a point in 
$\cals M^{\times}$, with $p$\ the
(single) intersection point. 
In suitable coordinates on $M$ $\dot{\si}(p)=\partial_1(p)$\ and
$\dot{\bar{\si}}(p)=\partial_d(p)$. 
We define local coordinates on $\cals M$ using 
$\psi(\cdot)=x^1(\cdot)-x^1(p)$\ and
$\bar{\psi}(\cdot)=x^d(\cdot)-x^d(p)$, as explained in Section
\ref{sec:strsp}. 
At $(\si,\bar{\si})$, we then find the local expression
\begin{equation}\label{eq:loccurr1}
  \begin{array}{rl}
  C^{\times}_{(\si,\bar{\si})}=&
  \sum_{k=1}^{d-1}\frac{(-1)^{k}}{(k-1)!(d-k-1)!} 
    \vep_{1\nu_2\ldots\nu_k\bar{\nu}_{k+1}\ldots\bar{\nu}_{d-1}d}\\
  & \td\si^{\nu_2}(0)\ldots\td\si^{\nu_k}(0)
    \td\bar{\si}^{\bar{\nu}_{k+1}}(0)\ldots
    \td\bar{\si}^{\bar{\nu}_{d-1}}(0)\\
  & \de(\si^2(0)-\bar{\si}^2(0))\ldots\de(\si^{d-1}(0)-\bar{\si}^{d-1}(0))
        \quad.
  \end{array}
\end{equation}
We must check that this is the current corresponding to $\cals M^{\times}$;
(see Appendix \ref{app:loc}).
\begin{itemize}
\item[(a)] $C^{\times}$ is localized on $\cals M^{\times}$, since, 
as one can see from (\ref{eq:loccurr}), it vanishes when the two strings do
not intersect. 
\item[(b)] A tangent vector, $v+\bar{v}$, at $(\si,\bar{\si})$ is parallel to
  $\cals M^{\times}$\ iff there exist real numbers $\al$\ and $\bar{\al}$\
  such that 
\begin{equation}\label{eq:transcond}
  v(0)+\al\dot{\si}(0)=\bar{v}(0)+\bar{\al}\dot{\bar{\si}}(0)
  \quad.
\end{equation}
A simple calculation shows that $C^{\times}$\ is transverse to $\cals
M^{\times}$, i.e., for all vectors $\pi=v+\bar{v}$\ fulfilling (\ref{eq:transcond}), one has
\begin{equation}\label{eq:iotaptheta}
  \iota_{\pi}C^{\times}_{(\si,\bar{\si})}=0 \quad.
\end{equation}
\item[(c)] Comparing (\ref{eq:loccurr1}) to equation (\ref{eq:formdec}) in
  Appendix \ref{app:loc}, we see that the regular part of $C^{\times}$ 
at $(\si,\bar{\si})$ is given by
\begin{equation}
  \begin{array}{rl}
  \widehat{C^{\times}_{(\si,\bar{\si})}}=&
  \sum_{k=1}^{d-1}\frac{(-1)^{k}}{(k-1)!(d-k-1)!} 
    \vep_{1\nu_2\ldots\nu_k\bar{\nu}_{k+1}\ldots\bar{\nu}_{d-1}d}\\
  & \td\si^{\nu_2}(0)\ldots\td\si^{\nu_k}(0)
    \td\bar{\si}^{\bar{\nu}_{k+1}}(\bar{0})\ldots
    \td\bar{\si}^{\bar{\nu}_{d-1}}(0)
        \quad,
  \end{array}
\end{equation}
and the localization functions are given by
\begin{equation}
  f_1=\si^2(0)-\bar{\si}^2(0)\quad \ldots\quad 
  f_{d-2}=\si^{d-1}(0)-\bar{\si}^{d-1}(0)
  \quad.
\end{equation}
It is easy to see that at $(\si,\bar{\si})$
\begin{equation}
  \left|\big{<}\;\cdot\;;\widehat{C^{\times}_{(\si,\bar{\si})}}\big{>}
        \right|=
  \left|\big{<}\;\cdot\;;
  \td(\si^2(0)-\bar{\si}^2(0))\ldots 
  \td(\si^{d-1}(0)-\bar{\si}^{d-1}(0))
  \big{>} 
        \right| 
  \quad.
\end{equation}
\end{itemize}

\begin{figure}
\begin{center}
\begin{picture}(160,60)
\put(22,30){\oval(24,20)[l]}
\put(28,30){\oval(24,20)[r]}
\put(40,17){\oval(30,7)[b]}
\put(40,43){\oval(30,7)[t]}
\put(25,23){\line(0,1){20}}
\put(55,17){\line(0,1){26}}
\put(26,25){\vector(0,1){10}}
\put(56,25){\vector(0,1){10}}

\put(25,23){\line(0,-1){6}}
\put(22,20){\line(1,0){6}}
\put(25,20){\circle*{1}}

\put(7,30){\makebox(0,0){$\si$}}
\put(58,30){\makebox(0,0){$\bar{\si}$}}

\put(77.5,30){\makebox(0,0){$\lmt$}}

\put(112,30){\oval(24,20)[l]}
\put(118,30){\oval(24,20)[r]}
\put(130,17){\oval(30,7)[b]}
\put(130,43){\oval(30,7)[t]}
\put(115,23){\line(0,1){20}}
\put(145,17){\line(0,1){26}}
\put(116,25){\vector(0,1){5}}
\put(146,25){\vector(0,1){5}}

\put(112,17){\oval(6,6)[tr]}
\put(118,23){\oval(6,6)[bl]}

\put(120,52){\makebox(0,0){$\PHI(\si,\bar{\si})$}}

\end{picture}

  \parbox{14cm}{
    \caption{\label{fig:concatenation}The map $\PHI$.}}
\end{center}
\end{figure}
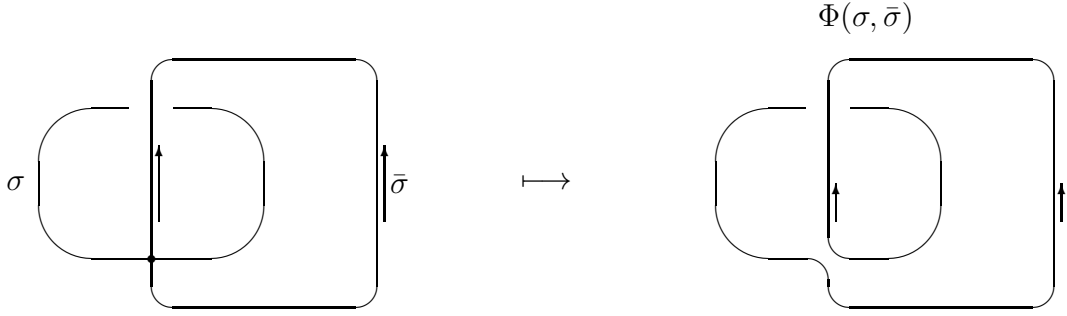

Let 
\begin{equation}\label{eq:PHI}
  \PHI:\cals M^{\times}\lra\cals M 
\end{equation}
be the map that associates to two intersecting strings 
their concatenation, with an appropriate scaling of the velocity
vectors, as shown in Figure \ref{fig:concatenation}. 
This map is nearly everywhere well-defined, namely on pairs of strings with one
self-intersection, but $n$-valued when the two strings intersect $n$\ times.

\begin{figure}[p!]
\begin{center}
  \input{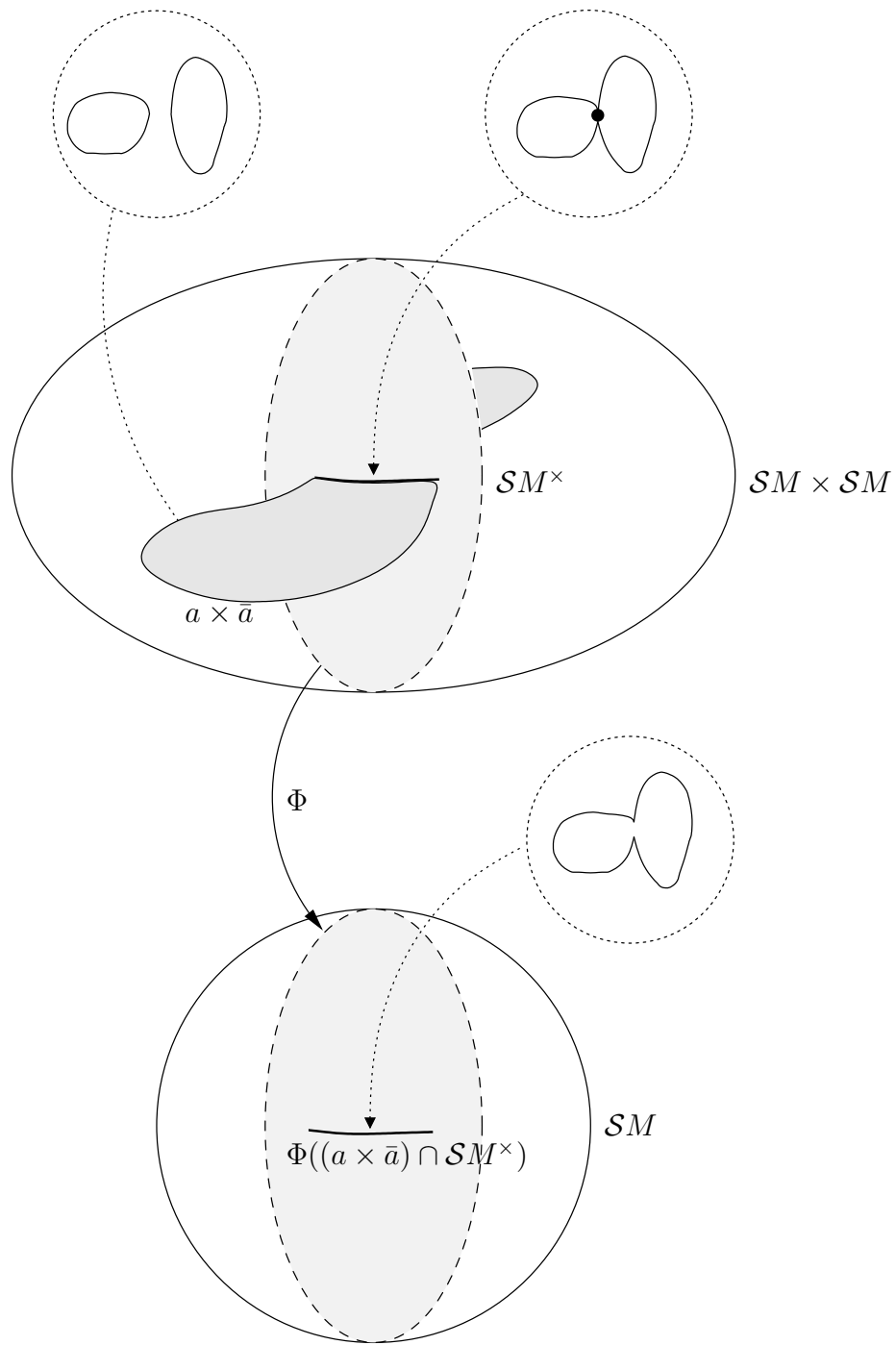}
  \parbox{14cm}{
    \caption{\label{fig:stringbracket}The definition of the string bracket.}}
\end{center}
\end{figure}

The string bracket is defined on string homology by
\footnote{Our definition differs from that described by Chas and Sullivan by a sign given by
$\{a;\bar{a}\}=\{\bar{a};a\}_{\mathrm{Chas-Sullivan}}$.}
(see also Figure \ref{fig:stringbracket})
\begin{equation}\label{eq:strbr}
  \begin{array}{cccc}
  \{\cdot;\cdot\}:& 
      \calh_{i}M\times\calh_{\bar{i}}M&\lra&\calh_{i+\bar{i}+2-d}M\\
      & (a,\bar{a}) & \lmt 
      & \{a;\bar{a}\}=
        (-1)^{\bar{i}(d+i)}
        \PHI\left((a\times\bar{a})\cap_{C^{\times}}\cals M^{\times}\right)
  \end{array}
  \quad.
\end{equation}
The r\^{o}le of $C^{\times}$\ is to orient the cycle obtained by intersecting 
an appropriately transversal representative $a\times\bar{a}$\ with 
$\cals M^{\times}$; see Appendix \ref{app:loc}.
The sign factor appearing in (\ref{eq:strbr}) is chosen in such a way that the
bracket is even, for even $d$, and odd, for odd $d$; in fact, it then
satisfies:  
\begin{itemize}
\item[(1)] \textit{Antisymmetry}
\begin{equation}
  \{a;\bar{a}\}=
    -(-1)^{(|a|+d)(|\bar{a}|+d)}\{\bar{a};a\}
    \quad,
\end{equation}
as can be checked by exchanging the factors in (\ref{eq:strbr}), and using 
$\mathrm{Ex}^{\ast}C^{\times}=(-1)^{1+d}C^{\times}$, with $\mathrm{Ex}$ 
the map that permutes the factors in $\cals M\times\cals M$.
\item[(2)]\textit{Jacobi identity}
\begin{equation}
  \{a;\{b;c\}\}=\{\{a;b\};c\}+(-1)^{(|a|+d)(|b|+d)}\{b;\{a;c\}\}
  \quad,
\end{equation}
(see Appendix \ref{app:jacobi} for a proof).
\end{itemize}
Here the degree $|\cdot|$\ of a cycle is its dimension.

Consider the symmetric algebra $S(\calh_{\ast}M)$\ over $\calh_{\ast}M$, with 
the grading given by $|\cdot|$. Extending the bracket as a superderivation, 
namely in such a way that the
\begin{itemize}
\item[(3)]\textit{Leibnitz rule}
\begin{equation}
  \{a,bc\}=\{a,b\}c+(-1)^{|b|(|a|+d)}b\{a,c\}  
\end{equation}
\end{itemize}
is fulfilled, one finds that $S(\calh_{\ast}M)$\ is \\
$\Diamond$\ a super-Poisson algebra (even bracket), for $d$\ even,\\ 
$\Diamond$\ a Gerstenhaber algebra (odd bracket), for $d$\ odd.

\section{A peculiarity of $\GLN$}
  \label{sec:gln}

In this section we highlight a property of $\GLN$\ which will be needed in
Section \ref{sec:homo}.

\vspace{4mm}
\noindent 
Let $G=\GLN$, and let $\rho$ denote its
standard representation.  
We define an invariant bilinear form $\ka$ as the trace in this representation:
\begin{equation}\label{eq:TT}
  \ka_{ab}=\ka(T_a,T_b)=\mathrm{tr}\left[\rho(T_a)\rho(T_b)\right]
  \quad.
\end{equation}
It then follows that
\begin{equation}
  \label{eq:TTexch0}
  \left(\ka^{ab}\rho(T_a)\otimes\rho(T_b)\right)v\otimes w= 
  w\otimes v
  \quad,
\end{equation}
where $v$\ and $w$\ are vectors in the representation space of $\rho$. In components with respect to a basis in this space the above identity reads
\begin{equation}
  \label{eq:TTexch}
  \left(\ka^{ab}\rho(T_a)^r_p\otimes\rho(T_b)^s_q\right)= 
  \de^r_q\de^p_s
  \quad.
\end{equation}
To prove this identity, we define a basis $\{E_{ij}|i,j=1..n\}$\ of $\gln$\ by 
setting $\rho(E_{ij})^r_s=\de_i^r\de_{js}$. For this basis, one finds that
$\ka(E_{ij},E_{kl})=\de_{il}\de_{jk}$. Equation (\ref{eq:TTexch}) then follows immediately.

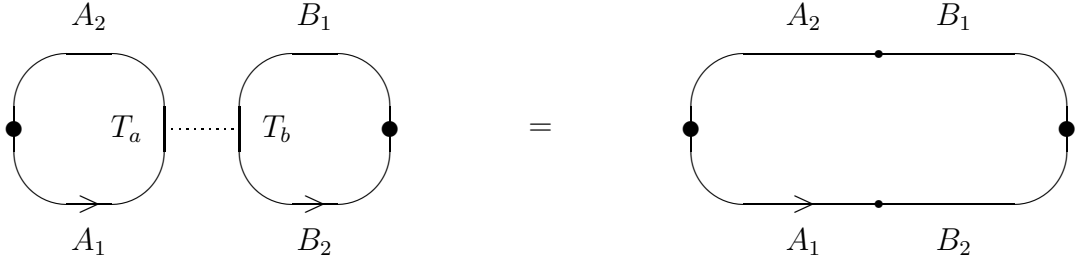
\begin{figure} 
\begin{picture}(160,50)(0,10)
\put(20,25){\oval(20,20)}
\put(50,25){\oval(20,20)}
\put(125,25){\oval(50,20)}
\qbezier[10](30,25)(35,25)(40,25)
\put(10,25){\circle*{2}}
\put(60,25){\circle*{2}}
\put(100,25){\circle*{2}}
\put(150,25){\circle*{2}}
\put(125,15){\circle*{1}}
\put(125,35){\circle*{1}}
\put(20,15){\makebox(0,0){$>$}}
\put(50,15){\makebox(0,0){$>$}}
\put(115,15){\makebox(0,0){$>$}}
\put(20,10){\makebox(0,0){$A_1$}}
\put(20,40){\makebox(0,0){$A_2$}}
\put(50,10){\makebox(0,0){$B_2$}}
\put(50,40){\makebox(0,0){$B_1$}}
\put(115,10){\makebox(0,0){$A_1$}}
\put(115,40){\makebox(0,0){$A_2$}}
\put(135,10){\makebox(0,0){$B_2$}}
\put(135,40){\makebox(0,0){$B_1$}}
\put(25,25){\makebox(0,0){$T_a$}}
\put(45,25){\makebox(0,0){$T_b$}}
\put(80,25){\makebox(0,0){$=$}}
\end{picture}

\begin{center}
  \parbox{14cm}{
    \caption{\label{fig:gln}Pictorial representation of (\ref{eq:trgln}).}}
\end{center}
\end{figure}

In the following, expressions of the form 
\begin{equation}\label{eq:HH}
  \mathrm{tr}_{\rho}\left[A_1T_aA_2\right]\ka^{ab}
  \mathrm{tr}_{\rho}\left[B_1T_bB_2\right]
\end{equation}
will appear, where $\rho$\ is a representation of $G$, $\{T_a\}$ is a 
basis of $\frag$, and 
$A_.,B_.$\ are elements of $\ssuu\frag$. 
For $G=\GLN$\ and $\rho$\ the standard representation,
such expressions can be simplified using (\ref{eq:TTexch}), as pictorially  
represented in Figure \ref{fig:gln}:
\begin{equation}\label{eq:trgln}
  \mathrm{tr}_{\rho}\left[A_1T_aA_2\right]\ka^{ab}
  \mathrm{tr}_{\rho}\left[B_1T_bB_2\right]=
  \mathrm{tr}_{\rho}\left[A_1B_2B_1A_2\right]
  \quad.
\end{equation}

\section{An algebra homomorphism 
        from $S(\calh_{\ast}M)$ to $H^{\ast}_{\de}\calo$}
  \label{sec:homo}

In this section we show that the map 
\begin{equation}\label{eq:homo0}
  \begin{array}{rrcl}
    S\ssh:
    & S(\calh_{\ast}M) & \lmt & H^{\ast}_{\de}\calo \\
    & a_1\ldots a_k & \lmt & 
      \big{<}a_1,\ssh_{\rho;A}\big{>}\ldots
      \big{<}a_k,\ssh_{\rho;A}\big{>}
  \end{array}
  \quad,
\end{equation}
which associates to a cycle in string
homology the corresponding observable of the \tft, based on the 
group $\GLN$\ in 
the standard representation, is a super-Poisson/Gerstenhaber algebra
homomorphism. This is accomplished by establishing the following properties:
\begin{eqnarray}
  i)&&
  |a|=\left|\big{<}a,\ssh_{\rho;A}\big{>}\right|
  \quad,\\
  \label{eq:homo}
  ii)&&
  \big{<}\{a;\bar{a}\},\ssh_{\rho;A}\big{>}=
  \{\big{<}a,\ssh_{\rho;A}\big{>};\big{<}\bar{a},\ssh_{\rho;A}\big{>}\}
  \quad.
\end{eqnarray}
Property $i)$\ follows from (\ref{eq:deghol}). 
Property $ii)$, is proven in several steps:

\noindent\textbf{Step 1}\\
Applying (\ref{eq:contraction1}), one finds that
\begin{equation}\label{eq:homo2}
  \{\big{<}a,\ssh\big{>};
       \big{<}\bar{a},\bar{\ssh}\big{>}\}=
  (-1)^{|\bar{a}|(d+|a|)}
  \big{<}a\times\bar{a},\{\ssh;\bar{\ssh}\}\big{>}
  \quad.
\end{equation}

\noindent\textbf{Step 2}\\
We derive a local expression for $\{\ssh,\bar{\ssh}\}$\ on 
$\ssll U\times\ssll U$.
First one verifies that
\begin{equation}
  \begin{array}{rl}
  & \left.\frac{d}{dt}\right|_{t=0} \ssh(C+t\eta) = \\
  =& \sum_{k=0}^d\int_{s=0}^{s=1}\mathrm{tr}\left[
    \left.\sshol(C)\right|_0^s \;
    \frac{1}{(k-1)!}\dot{\ga}^{\mu_1}(s)ds\td\ga^{\mu_2}(s)
    \ldots\td\ga^{\mu_k}(s)\eta_{\mu_1\mu_2\ldots\mu_k}(\ga(s))
    \; \left.\sshol(C)\right|_s^1\right]\quad.
  \end{array}  
\end{equation}
Using (\ref{eq:defder}), one finds the local expressions for
$\ddcr\ssh,\;\ssh\ddcl\in\calo(\OM^{\ast}(\ssll M)_{\cale}\cdot
M_{\cale}\cdot\ssuu\frag_{\cale})$, namely 
\begin{eqnarray}
  \ddcr\ssh & = & 
  \sum_{k=1}^{d}\frac{(-1)^{(k+1)(d+1)}}{(k-1)!(d-k)!}\int_{s=0}^{s=1}ds
  \; \de^{(d)}(\ga(s)-x)
  \vep_{\nu_1\nu_2\ldots\nu_k\mu_{k+1}\ldots\mu_d}\nonumber\\
  &  & 
  dx^{\mu_{k+1}}\ldots dx^{\mu_d}\dot{\ga}^{\nu_1}(s)
  \td\ga^{\nu_2}(s)\ldots\td\ga^{\nu_k}(s)\nonumber\\
  &  & 
  \mathrm{tr}
  \left[\left.\sshol(C)\right|_0^sT_a\left.\sshol(C)\right|_s^1\right]
  \otimes \ka^{ab}T_b
  \quad,
\end{eqnarray}
and
\begin{eqnarray}
  \ssh\ddcl& = & 
  \sum_{k=0}^{d-1}\frac{(-1)}{(k)!(d-k-1)!}\int_{s=0}^{s=1}ds
  \; \de^{(d)}(\ga(s)-x)
  \vep_{\mu_1\ldots\mu_k\nu_{k+1}\nu_{k+2}\ldots\nu_d}\nonumber\\
  &  & 
  dx^{\mu_1}\ldots dx^{\mu_k}
  \dot{\ga}^{\nu_{k+1}}(s)\td\ga^{\nu_{k+2}}(s)\ldots\td\ga^{\nu_d}(s)\nonumber\\
  &  & 
  \mathrm{tr}
  \left[\left.\sshol(C)\right|_0^sT_a\left.\sshol(C)\right|_s^1\right]
  \otimes \ka^{ab}T_b
  \quad.
\end{eqnarray}
Equation (\ref{eq:fieldbr}) yields the local expression for 
$\{\ssh;\bar{\ssh}\}\in\calo(\OM^{\ast}(\ssll M\times \ssll M)_{\cale})$
\begin{eqnarray} \label{eq:calcolohh}
  \{\ssh;\bar{\ssh}\} & = &
    \sum_{k=1}^{d-1}\frac{(-1)^{d+1}}{(k-1)!(d-k-1)!}
    \int_{s=0}^{s=1}ds\int_{\bar{s}=0}^{\bar{s}=1}d\bar{s}
    \; \de^{(d)}(\ga(s)-\bar{\ga}(\bar{s}))
    \; \vep_{\nu_1\nu_2\ldots\nu_k
      \bar{\nu}_{k+1}\bar{\nu}_{k+2}\ldots\bar{\nu}_d}
    \nonumber\\
    & &
    \dot{\ga}^{\nu_1}(s)\td\ga^{\nu_2}(s)\ldots\td\ga^{\nu_k}(s)
    \dot{\bar{\ga}}^{\bar{\nu}_{k+1}}(\bar{s})
    \td\bar{\ga}^{\bar{\nu}_{k+2}}(\bar{s})
    \ldots\td\bar{\ga}^{\bar{\nu}_d}(\bar{s}) \nonumber\\
    &&
    \mathrm{tr}\left[
      \left.\sshol\right|_0^sT_a\left.\sshol\right|_s^1\right]
    \ka^{ab}
    \mathrm{tr}\left[
      \left.\bar{\sshol}\right|_0^{\bar{s}}T_b
      \left.\bar{\sshol}\right|_{\bar{s}}^1\right]
    \quad.
\end{eqnarray}
We see that the latter can be written using the current $C^{\times}$, i.e.
\begin{equation}
  \{\ssh;\bar{\ssh}\} = C^{\times}\cdot
    \mathrm{tr}\left[
      \left.\sshol\right|_0^sT_a\left.\sshol\right|_s^1\right]
    \ka^{ab}
    \mathrm{tr}\left[
      \left.\bar{\sshol}\right|_0^{\bar{s}}T_b
      \left.\bar{\sshol}\right|_{\bar{s}}^1\right]
    \quad,
\end{equation}
which, for $G=\GLN$ in the standard representation, is equal to
\begin{equation}\label{eq:holcurrexp}
  \{\ssh;\bar{\ssh}\} = C^{\times}\cdot\sshh
  \quad,
\end{equation}
where
\begin{equation}
  \sshh=
   \mathrm{tr}\left[
    \left.\sshol\right|_0^s
    \left.\bar{\sshol}\right|_{\bar{s}}^1
    \left.\bar{\sshol}\right|_0^{\bar{s}}
    \left.\sshol\right|_s^1\right]
    \quad;
 \end{equation}
see (\ref{eq:trgln}). 

\noindent\textbf{Step 3}\\
$\{\ssh;\bar{\ssh}\}$ defines a form on $\cals M\times\cals M$.
{}From (\ref{eq:holcurrexp}) and (\ref{eq:ah}) one finds that
\begin{equation}
  \big{<}a\times\bar{a},\{\ssh;\bar{\ssh}\}\big{>}=
  \big{<}a\times\bar{a},C^{\times}\cdot\sshh\big{>}=
  \big{<}(a\times\bar{a})\cap_{C^{\times}}\cals M^{\times},\sshh\big{>}
  \quad.
\end{equation}
Moreover, one has that
\begin{equation}\label{eq:homo1}
  \big{<}\{a;\bar{a}\},\ssh\big{>}=
  (-1)^{|\bar{a}|(d+|a|)}
  \big{<}\PHI\left((a\times\bar{a})\cap_{C^{\times}}\cals M^{\times}\right),
         \ssh\big{>}
  \quad.
\end{equation}
Thus, to prove (\ref{eq:homo}), we simply have to show that
\begin{equation}
  \big{<}(a\times\bar{a})\cap_{C^{\times}}\cals M^{\times},\sshh\big{>}=
  \big{<}\PHI\left((a\times\bar{a})\cap_{C^{\times}}\cals M^{\times}\right),
         \ssh\big{>}
  \quad,
\end{equation}
which holds, as described in (\ref{eq:ph}), if
\begin{equation}
  \big{<}\Pi,
  \mathrm{tr}
  \left[
    \left.\sshol\right|_0^s
    \left.\bar{\sshol}\right|_{\bar{s}}^1
    \left.\bar{\sshol}\right|_0^{\bar{s}}
    \left.\sshol\right|_s^1\right]\big{>}_{(\si,\bar{\si})}=
  \big{<}\PHI_{\ast}\Pi,\ssh\big{>}_{\PHI(\si,\bar{\si})}
\end{equation}
for any $(\si,\bar{\si})\in\cals M^{\times}$\ and any parallel multivector
$\Pi\in\LA_{\ast}T_{(\si,\bar{\si})}\cals M^{\times}$. 
The validity of the latter follows immediately from the reparametrization
invariance of $\sshol$. 
The theorem is thus proven.

\section{Outlook}
  \label{sec:outlook}
In this section we outline various extensions and generalizations of the results
proven in this paper.
\subsection{Generalizations to other groups}
We start by describing some ideas about how to generalize the results of
this article by replacing $\GLN$\ with an arbitrary Lie group. 
Inspiration is taken from \cite{amr1}.\\

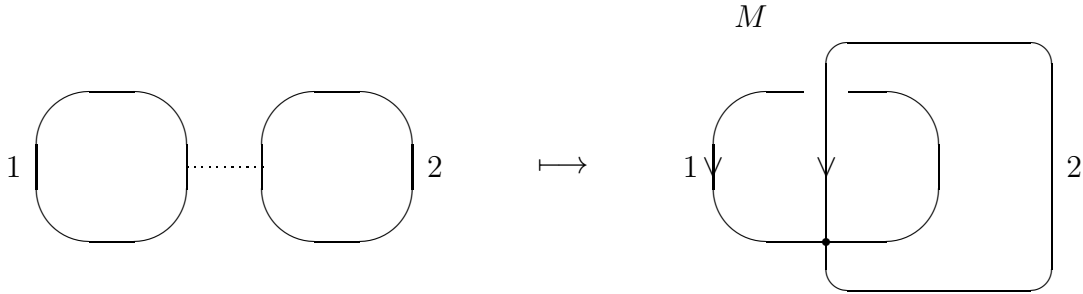
\begin{figure}[h]
\begin{center}
\begin{picture}(155,50)(0,5)
\put(20,25){\oval(20,20)}
\put(50,25){\oval(20,20)}
\qbezier[10](30,25)(35,25)(40,25)

\put(7,25){\makebox(0,0){$1$}}
\put(63,25){\makebox(0,0){$2$}}

\put(80,25){\makebox(0,0){$\lmt$}}

\put(112,25){\oval(24,20)[l]}
\put(118,25){\oval(24,20)[r]}
\put(130,12){\oval(30,7)[b]}
\put(130,38){\oval(30,7)[t]}
\put(115,18){\line(0,1){20}}
\put(145,12){\line(0,1){26}}
\put(100,25){\makebox(0,0){$\vee$}}
\put(115,25){\makebox(0,0){$\vee$}}

\put(115,18){\line(0,-1){6}}
\put(112,15){\line(1,0){6}}
\put(115,15){\circle*{1}}

\put(97,25){\makebox(0,0){$1$}}
\put(148,25){\makebox(0,0){$2$}}
\put(105,45){\makebox(0,0){$M$}}

\end{picture}

  \parbox{14cm}{
    \caption{\label{fig:chords}A chord diagram on $M$.}}
\end{center}
\end{figure}

\vspace{2mm}
A chord diagram (see Figure \ref{fig:chords}) 
is a union of disjoint oriented $S^1$-circles and disjoint arcs,
with the endpoints of the arcs on the circles.
A chord diagram on a manifold $M$\ (see Figure \ref{fig:chords})
is a (continuous) map from a chord diagram
to $M$\ such that each arc is mapped to a single point in $M$ (that is, 
each arc is mapped to an intersection of strings in $M$), modulo the obvious
action of $S^1$ on any circle. 
Let $\mathrm{ch}(M)$\ be the space of chord diagrams on $M$. 
It can be viewed
as a ``manifold'' with singularities when a circle is mapped to a single
point (just like for $\cals M$), and 
boundaries when two different crossings between circles approach one another
along one of the circles (see Figure \ref{fig:chordboundary}).

\begin{figure}[h]
\begin{center}
\begin{picture}(80,40)(-10,0)
\put(20,0){\line(0,1){40}}
\put(40,0){\line(0,1){40}}
\put(0,20){\line(1,0){60}}
\put(20,20){\circle*{1}}
\put(40,20){\circle*{1}}
\put(22,22){\vector(1,0){16}}
\end{picture}

  \parbox{14cm}{
    \caption{\label{fig:chordboundary}
      Approaching a boundary on $\mathrm{ch}(M)$ .}}
\end{center}
\end{figure}

One then defines a boundary operator, $\partial^{\mathrm{ch}(M)}$\, 
on cells in $\mathrm{ch}(M)$\ in such a way that
the so called $4T$-relation, represented in
Figure \ref{fig:fourt}, is respected.
\begin{figure}
\begin{center}
\begin{picture}(140,50)(0,10)
\put(0,30){\makebox(0,0){$\partial^{\mathrm{ch}}$}}
\put(10,15){\line(0,1){30}}
\put(10,15){\line(1,0){5}}
\put(10,45){\line(1,0){5}}
\put(135,15){\line(0,1){30}}
\put(135,15){\line(-1,0){5}}
\put(135,45){\line(-1,0){5}}
\put(145,30){\makebox(0,0){$=0$}}

\put(15,30){\makebox(0,0){$+$}}
\put(45,30){\makebox(0,0){$-$}}
\put(75,30){\makebox(0,0){$+$}}
\put(105,30){\makebox(0,0){$-$}}

\put(20,20){\vector(1,1){20}}
\put(40,20){\vector(-1,1){20}}
\put(25,20){\vector(0,1){20}}
\put(30,30){\circle*{1}}
\put(25,25){\circle*{1}}
\put(27,23){\vector(1,1){4}}

\put(50,20){\vector(1,1){20}}
\put(70,20){\vector(-1,1){20}}
\put(65,20){\vector(0,1){20}}
\put(60,30){\circle*{1}}
\put(65,35){\circle*{1}}
\put(63,37){\vector(-1,-1){4}}

\put(80,20){\vector(1,1){20}}
\put(100,20){\vector(-1,1){20}}
\put(95,20){\vector(0,1){20}}
\put(90,30){\circle*{1}}
\put(95,25){\circle*{1}}
\put(93,23){\vector(-1,1){4}}

\put(110,20){\vector(1,1){20}}
\put(130,20){\vector(-1,1){20}}
\put(115,20){\vector(0,1){20}}
\put(120,30){\circle*{1}}
\put(115,35){\circle*{1}}
\put(117,37){\vector(1,-1){4}}

\end{picture}
  \parbox{14cm}{
    \caption{\label{fig:fourt}$4T$-relations.}}
\end{center}
\end{figure}

The chord homology $\calh_{\ast}^{\mathrm{ch}}M$\ is the homology of
$\mathrm{ch}(M)$\ with respect to $\partial^{\mathrm{ch}(M)}$.

In analogy to $\cals M^{\times}\subset\cals M\times\cals M$\, one defines 
$\mathrm{ch}(M)^{\times}\subset\mathrm{ch}(M)\times\mathrm{ch}(M)$\ as
the space of pairs of chord diagrams on $M$\ whose strings intersect at least
once. 
Similarly to $\PHI:\cals M^{\times}\lra\cals M$\, one defines the (generally
multivalued) map
\begin{equation}
  \PHI^{\mathrm{ch}}:\mathrm{ch}(M)^{\times}\lra\mathrm{ch}(M)
\quad,  
\end{equation}
which associates to a pair of chord diagrams on $M$\ with one 
intersection point the 
union of the two chord diagrams with a new arc corresponding to the
intersection (and in an analogous way for multiple intersection points).

As in equation (\ref{eq:strbr}), one defines a bracket
\begin{equation}\label{eq:chordbr}
  \begin{array}{cccc}
  \{\cdot;\cdot\}:& 
      \calh^{\mathrm{ch}}_{i}M\times\calh^{\mathrm{ch}}_{\bar{i}}M
      &\lra&\calh^{\mathrm{ch}}_{i+\bar{i}+2-d}M \\
      & (a,\bar{a}) & \lmt 
      & (-1)^{\bar{i}(i+d)}\PHI^{\mathrm{ch}}
        ((a\times\bar{a})\cap_{C^{\times}}\mathrm{ch}(M)^{\times})
  \end{array}
  \quad,
\end{equation}
which is a bracket/antibracket for $d$\ even/odd; the current $C^{\times}$\ on $\mathrm{ch}(M)^{\times}$\ can be constructed in a similar way as in section \ref{sec:bracket}.

Similarly to $S(\calh_{\ast}M)$, it is possible to define a
super-Poisson/Gerstenhaber algebra $S(\calh_{\ast}^{\mathrm{ch}}M)$.

\begin{figure}
\begin{center}

\begin{picture}(140,60)
\put(10,10){\vector(1,0){40}}
\put(10,50){\vector(1,0){40}}
\qbezier[40](30,10)(30,30)(30,50)
\put(30,53){\makebox(0,0){$T_a$}}
\put(30,7){\makebox(0,0){$T_b$}}
\put(33,30){\makebox(0,0){$\ka^{ab}$}}
\put(27,47){\makebox(0,0){$t_1$}}
\put(27,13){\makebox(0,0){$t_2$}}
\put(7,50){\makebox(0,0){$\rho_1$}}
\put(7,10){\makebox(0,0){$\rho_2$}}

\put(65,30){\makebox(0,0){$\leadsto$}}
\put(110,50){\makebox(0,0){$
    \mathrm{tr}_{\rho_1}\left[\ldots
      \left.\sshol\right|_{\cdot}^{t_1}T_a\left.\sshol\right|^{\cdot}_{t_1}
      \ldots\right]$}}
\put(110,30){\makebox(0,0){$\ka^{ab}$}}
\put(110,10){\makebox(0,0){$
    \mathrm{tr}_{\rho_2}\left[\ldots
      \left.\sshol\right|_{\cdot}^{t_2}T_b\left.\sshol\right|^{\cdot}_{t_2}
      \ldots\right]$}}

\end{picture}

  \parbox{14cm}{
    \caption{\label{fig:ffcd}The form $\ssh$\ associated to (a part of) 
      a chord diagram; $t_1$ and $t_2$ refer to a $S^1$-parametrization of the
      circles.}} 
\end{center}
\end{figure}
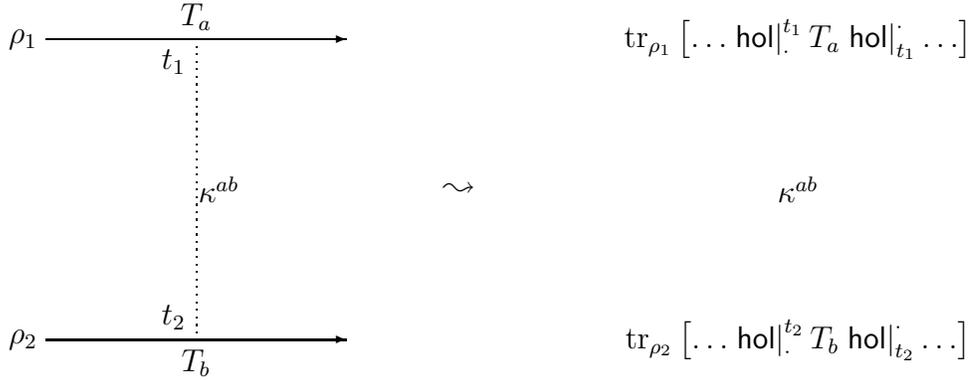

In analogy to (\ref{eq:homo0}), we define a map
\begin{equation}
  \begin{array}{rrcl}\label{eq:homoch}
    S(\ssh^{\mathrm{ch},G}):
    &S(\calh^{\mathrm{ch}}_{\ast}M^G) & \lmt & H^{\ast}_{\de}\calo^G \\
    & a_1\ldots a_k & \lmt & 
      \big{<}a_1,\ssh^{\mathrm{ch}}\big{>}\ldots
      \big{<}a_k,\ssh^{\mathrm{ch}}\big{>}
  \end{array}
  \quad,
\end{equation}
where $\calh^{\mathrm{ch}}_{\ast}M^G$\ denotes the homology of chord diagrams
with circles labeled by representations of $G$. The form
$\ssh^{\mathrm{ch}}$\ is defined as explained in Figure \ref{fig:ffcd}.

The map (\ref{eq:homoch}) is a super-Poisson-/Gerstenhaber 
algebra homomorphism. 
This can be proved by the same reasoning as that in Section 
\ref{sec:homo} and in \cite{amr1}.

The content of \cite{amr1} concerns the special case of the above construction for 
manifolds $M$\ of dimension $d=2$\ and for
$\calh^{\mathrm{ch}}_{0}M\subset\calh^{\mathrm{ch}}_{\ast}M$.

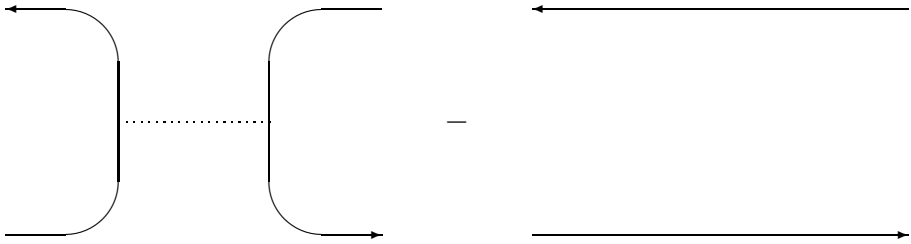
\begin{figure}[h]
\begin{center}
\begin{picture}(140,40)
\put(10,25){\oval(30,30)[r]}
\put(60,25){\oval(30,30)[l]}
\put(15,40){\vector(-1,0){5}}
\put(55,10){\vector(1,0){5}}

\qbezier[20](25,25)(35,25)(45,25)

\put(70,25){\makebox(0,0){$-$}}

\put(80,10){\vector(1,0){50}}
\put(130,40){\vector(-1,0){50}}

\end{picture}
  \parbox{14cm}{
    \caption{\label{fig:glnideal}The ``$\GLN$''-ideal $I$.}}
\end{center}
\end{figure}

The symmetric algebra on string homology, $S(\calh_{\ast}M)$,
  is obtained by taking the quotient of
$S(\calh^{\mathrm{ch}}_{\ast}M)$\ by the ideal $I$\ generated by the diagrams
of Figure \ref{fig:glnideal}.
\\

One then sees that the following diagram is commutative:
\begin{center}
\begin{picture}(120,40)
\put(30,30){\makebox(0,0){$S(\calh^{\mathrm{ch}}_{\ast}M)$}}
\put(30,10){\makebox(0,0){$S(\calh_{\ast}M)$}}
\put(30,26){\vector(0,-1){12}}
\put(33,20){\makebox(0,0){$\pi_I$}}
\put(40,30){\vector(1,0){30}}
\put(40,12){\vector(2,1){30}}

\put(55,33){\makebox(0,0){$S(\ssh^{\mathrm{ch},\GLN})$}}
\put(85,30){\makebox(0,0){$H^{\ast}_{\de}\calo^{\GLN}$}}

\end{picture}
\end{center}

\subsection{Generalization to nontrivial principal bundles}

In this subsection we explain how to extend methods and results of this paper
to the situation where $P$\ is a non-trivial bundle with base space $M$ and thus not necessarily admits a flat connection.

A principal bundle is determined by its ``transition functions''
\begin{equation}
  t_{ij}:U_i\cap U_j\lra G
\end{equation}
defined on intersections of two coordinate patches of $M$, and with the
property that
\begin{equation}
  t_{ij}t_{jk}=t_{ik}\quad\mathrm{on}\;U_i\cap U_j\cap U_k\quad.
\end{equation}
Two sets of transition functions $t,\tilde{t}$\ describe the same bundle iff
there exist ``gauge transformations'' 
\begin{equation}
  g_i:U_i\lra G
\end{equation}
such that
\begin{equation}
  t_{ij}=g_i\tilde{t}_{ij}g_j^{-1}\quad,\quad\mathrm{on}\;U_i\cap U_j\quad.
\end{equation}
A connection on $P$\ associates to every patch a $\frag$-valued one-form
\begin{equation}
  A_i\in\OM^1(U_i)\otimes\frag
  \quad,
\end{equation}
such that 
\begin{equation}
  A_i=t_{ij}A_jt_{ij}^{-1}+t_{ij}dt_{ij}^{-1}
      \quad\mathrm{on}\;U_i\cap U_j\quad.
\end{equation}
The curvature, $F$, of the connection $A$ is given, on every patch, by a
$\frag$-valued two-form 
\begin{equation}
  F_i=dA_i+\frac{1}{2}[A_i,A_i]\in\OM^2(U_i)\otimes\frag
  \quad,
\end{equation}
such that 
\begin{equation}
  F_i=t_{ij}F_jt_{ij}^{-1}\quad,\quad\mathrm{on}\;U_i\cap U_j\quad.
\end{equation}
The forms $C$ are $\frag_{\cale}$-valued forms. On every coordinate patch,
$C$ is given by
\begin{equation}
  C_i\in\OM^{\ast}(U_i)\otimes\frag_{\cale}
  \quad,
\end{equation}
with the property that 
\begin{equation}
  C_i=t_{ij}C_jt_{ij}^{-1}\quad,\quad\mathrm{on}\;U_i\cap U_j\quad.
\end{equation}
A principal bundle is trivial iff one can choose trivial transition functions:
$t_{ij}=\id$, for all $U_i,U_j$, with $U_i\cap U_j\neq\emptyset$.
The connection, the curvature and the forms $C$\ are then globally defined on
$M$. 

We now turn our attention to the master action and the bracket of the topological field theory.
The forms on the patches
\begin{equation}
  s_i=\mathrm{tr}_{\rho}
      \left[C_i(F_i+\frac{1}{2}d_{A_i}C_i+\frac{1}{3}C_i^2)\right]
      \in\OM^{\ast}(U_i)_{\cale}
\end{equation}
satisfy $s_i=s_j$\ on $U_i\cap U_j$, and thus yield a globally defined form
$s$\ on $M$.
We may therefore define a master action, $S$, by 
\begin{equation}
  S=\int_Ms=\int_M\mathrm{tr}_{\rho}
    \left[C(F+\frac{1}{2}d_AC+\frac{1}{3}C^2)\right]
      \in\OM^{\ast}(M)_{\cale}
  \quad.
\end{equation}
The bracket is well defined, since one has
\begin{equation}
  \{C_i;C_i\}=\{C_j;C_j\}
  \quad,
\end{equation}
a consequence of the invariance of the bilinear form $\ka$\ under the adjoint
action of $G$\ on $\frag$. 
The master action still satisfies the master equation $\{S;S\}=0$. Furthermore,
\begin{equation}\label{eq:deltafields1}
  \{S;C_i\}=\de C_i=(-1)^d (F_i+d_{A_i}C_i+C_i^2)
  \quad.
\end{equation}

We now address the task of defining  generalized parallel transporters and
generalized Wilson loops.
They can be defined as elements of $\LA^{\ast}T_{\ga}\ssll M$, for each loop
$\ga\in\ssll M$. 
Let $0=t_0<t_1<\ldots<t_{k-1}<t_k=1$, and let $U_1,\ldots,U_k=U_1$ be patches
on $M$ such that $\ga(t)\in U_i$, for $t\in[t_{i-1},t_i]$. One then defines
the trace the generalized Wilson loop as
\begin{equation}\label{eq:splithol}
  \ssh_{\rho;A}(C)=\mathrm{tr}_{\rho}
    \left[
    \left.\sshol_{A_1}(C_1)\right|_{0}^{t_1}
    t_{12}
    \left.\sshol_{A_2}(C_2)\right|_{t_1}^{t_2}
    \ldots
    \left.\sshol_{A_{k-1}}(C_{k-1})\right|_{t_{k-2}}^{t_{k-1}}
    t_{k-1,k}
    \left.\sshol_{A_k}(C_k)\right|_{t_{k-1}}^{1}
    \right]
  \quad.
\end{equation}
The factors $\left.\sshol_{A_i}(C_i)\right|_{t_{i-1}}^{t_i}$\ are defined as in (\ref{eq:parthol}). It is easy to see that this definition does not depend on the choice of the charts and is invariant under gauge transformations.
One then shows that
\begin{equation}\label{eq:dhol1}
  \td\ssh_{\rho;A}=\int_0^1d\tau \; \mathrm{tr}_{\rho}\left[
        \left.\sshol_{A}(C)\right|_{0}^{\tau}
        \iota_{\dot{\ga}}
        \mathrm{ev}^{\ast}_{\tau}(F+d_AC+C^2)
        \left.\sshol_{A}(C)\right|_{\tau}^{1} \right]
  \quad,
\end{equation}
where the $\tau$-integral has to be split, as in (\ref{eq:splithol}), if the loop crosses different patches.
Comparing (\ref{eq:splithol}) and (\ref{eq:deltafields1}), one finds that the fundamental identity (\ref{eq:ddeltah}) is fulfilled:
\begin{equation}
  ((-1)^d\de+\td)\ssh_{\rho;A}=0
  \quad.
\end{equation}

\subsection{Remarks on quantization}
The construction we have described in this paper yields, in the case of
an even-dimensional manifold $M$, a Poisson algebra
of observables (related to the string topology of $M$ if we choose
$GL(n)$ as our Lie group). It is then natural to ask if and how this Poisson
algebra may be quantized. We sketch in this Section a few approaches that
might help understanding this problem.

\subsubsection{Path-integral quantization}
If $d=\dim M$ is even, our approach describes the BRST formalism 
for a field theory in the Hamiltonian formalism with
the functional $S^{[d]}$ as the BRST  
generator.
If we want to quantize this theory using path-integrals,
we must first move to the Lagrangian formalism. As explained in 
Appendix~\ref{app:bvbrst},
the corresponding action functional on $N=M\times I$ is $S^{[d+1]}$.

In the case $d=2$, this is the BV action for Chern--Simons theory, and 
this is in accordance with the fact that Chern--Simons
theory provides a quantization of the Goldman \cite{gol} bracket (the
$2$-dimensional version of the string bracket), see \cite{amr}.
In higher dimensions,  $S^{[d+1]}$ defines new topological quantum field
theories (TQFT), 
among which we have the so-called $BF$ theories \cite{sch,blth}
which
can be obtained by
particular choices of the metric Lie algebra.

Our observables for strings on $M$ have then to be lifted to the corresponding
observables on $N=M\times I$ (or, more generally, on a 
$(d+1)$-dimensional manifold
$N$). The formulae we have given in odd dimensions describe this algebra of 
observables. Notice however that, in order to avoid singularities in the 
computation of expectation values, one has to restrict oneself
to imbedded strings in $N$ (and possibly also to introduce a framing).
In the particular case of $BF$ theories,
the expectation values of these observables correspond to the cohomology 
classes of
imbedded strings considered in \cite{ccl}, as shown in \cite{accr1,accr2}.
As a consequence, the quantization
of the string topology of $M$ must be related to the homology of the space
of imbedded strings in $M\times I$. 
This space must then be endowed with the structure
of associative algebra in such a way that its commutator yields, in the 
classical limit, the Poisson bracket of the projections of the strings to $M$.

\subsubsection{Deformation quantization}
For $d=2$ and $M$ non-compact,
the ideas described above have an explicit realization in
terms of deformation quantization (i.e., working with formal power
series in $\hbar$), as described in \cite{amr}. 
The construction is based on the Kontsevich integral for link
invariants \cite{konts-links} which is the perturbative formulation
of Chern--Simons theory in the holomorphic gauge studied in \cite{fk}.

The higher-dimensional generalization of this
approach should be obtained 
by considering perturbative expansions, in a suitable gauge,
of the corresponding TQFTs.

\subsubsection{Geometric quantization}
In some cases (e.g., $BF$ theories),
the Poisson subalgebra of functionals commuting with $S^{[d]}$ is
the algebra of a
reduced phase space of generalized gauge fields
on $M$. This space inherits a symplectic structure and one may try to quantize
it using deformation quantization and produce a TQFT in Atiyah's sense.
In the $2$-dimensional case, when the
reduced phase space turns out to be the space of flat connections on $M$
modulo gauge transformations, this program works 
(at least for compact groups).
One may regard quantum groups
as one of its outcomes. It would be very interesting to understand if the 
higher-dimensional case produces interesting generalizations thereof.


\begin{appendix}
\section{Intersection of cycles and currents}
\label{app:loc}

In this section we explain some concepts and manipulations used in the proof
of eq. (\ref{eq:homo}) in Section \ref{sec:homo}.

\vspace{4mm}
Let $A$\ be a manifold and $A^{\times}$\ an oriented immersion of codimension
$n$, which defines an element of the homology, $H_{\ast}A$, of $A$. 
Let $C^{\times}$\ be the current that localizes on this immersion, i.e., a
singular $n$-form on $A$ with the following properties: 
\begin{itemize}
\item[(a)] The form localizes on $A^{\times}$, i.e. for any point $p$\ not in $A^{\times}$\ one has
\begin{equation}
  C^{\times}_p=0
  \quad.
\end{equation}
\item[(b)] The form is transverse, i.e., for every point $p$\ on $A^{\times}$\
  and an arbitrary parallel tangent vector $P(p)\in T_pA^{\times}$\ one has
\begin{equation}
  \iota_{P(p)}C^{\times}_p=0
  \quad.
\end{equation}
\item[(c)] Let $A^{\times}$\ be defined, locally, as the zero-set of functions
  $f_1,\ldots,f_n$, with $df_1\ldots df_n\neq0$. Then the current $C^{\times}$
  is given by
\begin{equation} \label{eq:formdec}
  C^{\times}=\widehat{C^{\times}}\de(f_1)\ldots\de(f_n)
  \quad,
\end{equation}
where $\widehat{C^{\times}}$ is a regular form, and for every point $p$\ in $A^{\times}$\ and every multivector $V\in T_pA$,
\begin{equation}
  \left|\big{<}V,\widehat{C^{\times}}\big{>}\right|=
  \left|\big{<}V,df_1\ldots df_n\big{>}\right|
  \quad.
\end{equation}
\end{itemize}
In particular, $C^{\times}$\ defines an orientation on the normal bundle,
$N(A^{\times})$, of $A^{\times}$\ in $A$. 
Given an $i$-cycle $a\in H_{i}(A)$, one can define a new cycle by considering
the intersection 
\begin{equation}
  a\cap_{C^{\times}}A^{\times}\in H_{i-n}(A)
  \quad.
\end{equation}
As a set, it is obtained by intersecting an appropriate representative of $a$\
with $A^{\times}$. 
The orientation is defined as follows: Let $p$\ be a point in this intersection, $P\in\LA_{i-n}T_p(a\cap A^{\times})$\ the multivector that is the infinitesimal version at $p$\ of $a\cap A^{\times}$, $T\in\LA_{n}T_pa$\ the multivector in the normal bundle to $A^{\times}$\ such that $T\wedge P$\ is the infinitesimal version of $a$\ at $p$. Then one defines
\begin{equation}
  \ori_{a\cap_{C^{\times}}A^{\times}}(P)=
  \ori_a(T\wedge P)\cdot\ori_{N(A^{\times})}(T)
  \quad,
\end{equation}
where $\ori_{N(A^{\times})}$\ is given by the current $C^{\times}$.

For any closed form $H$ on $A$, one has that
\begin{equation}\label{eq:ah}
  \big{<}a,C^{\times}H\big{>}=
  \big{<}a\cap_{C^{\times}}A^{\times},H\big{>}
  \quad.
\end{equation}
Next, let $\PHI$\ be a map from $A^{\times}$\ into some other manifold $B$,
and $h$\ a closed form on $B$. If for an arbitrary point $p$\ in $A^{\times}$\
and any parallel multivector $P\in \LA_{\ast}T_pA^{\times}$, one has that
\begin{equation}\label{eq:ph}
  \big{<}P,H\big{>}_p=
  \big{<}\PHI_{\ast}P,h\big{>}_{\PHI(p)}
  \quad,
\end{equation}
then
\begin{equation}
  \big{<}a\cap_{C^{\times}}A^{\times},H\big{>}=
  \big{<}\PHI(a\cap_{C^{\times}}A^{\times}),h\big{>}
  \quad.
\end{equation}

\section{The Jacobi identity for the string bracket}
  \label{app:jacobi}

In this appendix we show how to prove the Jacobi identity for the string bracket of section \ref{sec:bracket}.

We first rewrite the Jacobi identity as
\begin{equation}\label{eq:appjac}
  (-1)^{\eta(abc)}\{\{a;b\};c\}+\mathrm{cycl.}(abc)=0
  \quad,
\end{equation}
where the sign factor is $\eta(abc)=(|a|+d)(|c|+d)$.
We can define the first term as
\begin{equation}\label{eq:appjac1}
  \begin{array}{c}
  (-1)^{\eta(abc)}\{\{a;b\};c\} = \\ \\
     =(-1)^{\eta(abc)+\si(abc)}\\ \left[
     \PHI^{(1,23)}\left(
     \left(a^{(1)}\times b^{(2)}\times c^{(3)}\right)
        \cap_{C^{\times(12)}\wedge C^{\times(13)}}
        \left(\cals M^{\times(12)}\cap\cals M^{\times(13)}\right)
     \right)+ \right. \\ \left.
     +\PHI^{(2,13)}\left(
     \left(a^{1}\times b^{2}\times c^{3}\right)
        \cap_{C^{\times(12)}\wedge C^{\times(23)}}
        \left(\cals M^{\times(12)}\cap\cals M^{\times(23)}\right)
     \right)
     \right]
  \quad.
  \end{array}
\end{equation}
Let us first explain the objects that appear in the above definition. The sign
factor is $\si(abc)=(|b|(d+|a|)+|c|(|a|+|b|)$, which follows from the
definition of the string bracket, (\ref{eq:strbr}). 
$a^{(1)}\times b^{(2)}\times c^{(3)}$\ is a cycle in $\cals M^{(1)}\times\cals
M^{(2)}\times\cals M^{(3)}$. 
A point in $\cals M^{\times(ij)}$\ is a triple of strings, $(\si_1,\si_2,\si_3)\in\cals M^{(1)}\times\cals M^{(2)}\times\cals M^{(3)}$, such that the $i$-th and the $j$-th intersect at least once. $C^{\times(ij)}$\ is the corresponding current. $\PHI^{(i,jk)}$\ is the map
\begin{equation}
  \PHI^{(i,jk)}:
  \cals M^{\times(ij)}\cap\cals M^{\times(ik)}\lra\cals M^{\times}
\end{equation}
which opens the intersections between the $i$-th and the $j$-th and between
the the $i$-th and the $k$-th string, in the same way as the map $\PHI$\
in (\ref{eq:PHI}) does.

Now consider the two terms appearing in (\ref{eq:appjac}) corresponding to the
cycle $a$\ intersecting both the cycles $b$\ and $c$. 
The first term  corresponds to the first term in (\ref{eq:appjac1}). 
The second one appears in $(-1)^{\eta(cab)}\{\{c;a\};b\}$\ and reads
\begin{equation}\label{eq:appjac2}
  \begin{array}{c}
     (-1)^{\eta(cab)+\si(cab)} \\
     \PHI^{(2,13)}\left(
     \left(c^{(1)}\times a^{(2)}\times b^{(3)}\right)
        \cap_{C^{\times(12)}\wedge C^{\times(23)}}
        \left(\cals M^{\times(12)}\cap\cals M^{\times(23)}\right)
     \right)
  \end{array}
  \quad.
\end{equation}
To prove that the Jacobi identity holds, we only have to prove that two such
terms add up to zero. 

We first write the second term, rearranging the indices and bringing the
cycles into a convenient order, i.e., 
\begin{equation}\label{eq:appjac3}
  \begin{array}{c}
     (-1)^{\eta(cab)+\si(cab)}(-1)^{|c|(|a|+|b|)}\\
     \PHI^{(1,32)}\left(
     \left(a^{(1)}\times b^{(2)}\times c^{(3)}\right)
        \cap_{C^{\times(31)}\wedge C^{\times(12)}}
        \left(\cals M^{\times(12)}\cap\cals M^{\times(23)}\right)
     \right)  \quad;
  \end{array}
\end{equation}
then we bring the currents into a convenient form
\begin{equation}\label{eq:appjac4}
  \begin{array}{c}
     (-1)^{\eta(cab)+\si(cab)}+(-1)^{|c|(|a|+|b|)+1}\\
     \PHI^{(1,23)}\left(
     \left(a^{(1)}\times b^{(2)}\times c^{(3)}\right)
        \cap_{C^{\times(12)}\wedge C^{\times(13)}}
        \left(\cals M^{\times(12)}\cap\cals M^{\times(23)}\right)
     \right) \quad,
  \end{array}
\end{equation}
using that $|C^{\times(ij)}|=d$\ and
$C^{\times(ij)}=(-1)^{d+1}C^{\times(ji)}$.
What remains to be shown is thus that
\begin{equation}
  \eta(abc)+\si(abc)+\eta(cab)+\si(cab)+|c|(|a|+|b|)+1\stackrel{!}{=}1
  \quad,
\end{equation}
which is easily seen to hold.

\section{Local expression for the generalized \\ parallel transporters}
 \label{app:lochol}

In local coordinates
$(\ga^{\mu}(t))_{t\in S^1}$\ the generalized holonomy reads
\begin{equation}
  \begin{array}{rcl}
  \left.\sshol^n_{A}(C)\right|_{t_i}^{t_f}& = & 
  \sum_{n_1,\ldots,n_n=1}^{\infty}
   \int_{(t_1,\ldots,t_n)\in\left.\DE_n\right|_{t_i}^{t_f}}\\
      &&\left.\mathrm{hol}_{A}\right|_{t_i}^{t_1}
      \dot{\ga}^{\mu^1_1}(t_1)dt_1
        \td\ga^{\mu^1_2}(t_1)\ldots\td\ga^{\mu^1_{n_1}}(t_1)
        C_{\mu^1_1\mu^1_2\ldots\mu^1_{n_1}}(\ga(t_1))
      \left.\mathrm{hol}_{A}\right|_{t_1}^{t_2}\\
      &&\ldots\\
      &&
      \left.\mathrm{hol}_{A}\right|_{t_{n-1}}^{t_n}
      \dot{\ga}^{\mu^n_1}(t_n)dt_n
        \td\ga^{\mu^n_2}(t_n)\ldots\td\ga^{\mu^n_{n_n}}(t_n)
        C_{\mu^n_1\mu^n_2\ldots\mu^n_{n_n}}(\ga(t_n))
      \left.\mathrm{hol}_{A}\right|_{t_n}^{t_f}
   \end{array}
  \quad,
\end{equation}
where $\td$\ is the differential on $\ssll M$.

\section{BV/BRST}
 \label{app:bvbrst}
In this Appendix we explain the relationship between $S^{[d+1]}$\ and
$S^{[d]}$, where $d$\ is an even number. We follow \cite{ht}.
For notational simplicity, we omit the Lie algebra part of the forms.

\vspace{4mm}
Let $N\ni\xx$\ be an oriented manifold with $\dim N=d$\ even,  and $M=I\times
N\ni(t,\xx)$\ with the product orientation. 
Let us write the fields on $M$\ as
\begin{equation}
  C=dtC_t+D=
    \sum_{k=0}^d dt\;\frac{1}{k!}\;dx^{i_1}\ldots dx^{i_k}C_{ti_1\ldots i_k}+
    \sum_{k=0}^d \frac{1}{k!}\;dx^{i_1}\ldots dx^{i_k}D_{i_1\ldots i_k}
    \quad.
\end{equation}
{}From (\ref{eq:brcc}) it follows that, in the BV-formalism,
 one can choose as fields and corresponding antifields, respectively,
\begin{equation}
  D_{i_1\ldots i_k} \quad\longleftrightarrow\quad 
  \frac{1}{(d-k)!}\;\vep^{i_1\ldots i_ki_{k+1}\ldots i_d}C_{ti_{k+1}\ldots i_d}
  \quad.
\end{equation}
After choosing a gauge in which the connection $A$\ has vanishing time component, $A_t=0$, the master action in the Lagrangian formalism reads
\begin{equation}
  S^{[d+1]}[C^t,D]=
    \int_I dt \int_N (d_A D+D^2)C_t+\frac{1}{2}\dot{D}D
  \quad.
\end{equation}
A gauge-fixing functional $\PSI[D]$\ ($|\PSI|=1$) defines a gauge-fixed action
\begin{equation}
  S^{[d+1]}_{\PSI}[D]=
             S^{[d+1]}
             \left[C_{ti_{k+1}\ldots i_d}=\frac{1}{k!}\;
             \vep^{i_1\ldots i_ki_{k+1}\ldots i_d}
             \frac{\stackrel{\ra}{\de}}{\de D_{i_1\ldots i_k}}\PSI,D\right]
  \quad.
\end{equation}
For a gauge-fixing functional adapted to the ``space-time'' split 
$M=I\times N$ of the form
\begin{equation}
  \PSI[D]=-\int_Idt \;K[D] \quad,
\end{equation}
where $K$\ is some functional of $D$, with $D$ interpreted as a form on $N$,
one finds that
\begin{equation}\label{eq:hamac}
  S^{[d+1]}_{\PSI}[D]=
      \int_I dt
      \left(\{S^{[d]},K\}_N+\frac{1}{2}
      \int_N \dot{D}D\right)
    \quad,
\end{equation}
with
\begin{equation}
  S^{[d]}[D]=\frac{1}{2}\int_N Dd D+\frac{2}{3}D^3
  \quad.
\end{equation}
We remark that the gauge fixed action (\ref{eq:hamac}) is already in 
Hamiltonian form, since it is of first order in time derivatives.
Since
\begin{equation}
  \left.\{S^{[d+1]};D_{i_1\ldots i_k}(t,\xx)\}\right|
  _{C_t=\frac{\stackrel{\ra}{\de}}{\de D}\PSI}
  =(-1)^k(dD+D^2)_{i_1\ldots i_k}(t,\xx)
\end{equation}
and
\begin{equation}
  \{S^{[d]};D_{i_1\ldots i_k}(\xx)\}=
        (-1)^k(dD+D^2)_{i_1\ldots i_k}(\xx)
  \quad,
\end{equation}
$S^d$\ can be interpreted as the BRST-generator in the Hamiltonian
formalism, and (\ref{eq:hamac}) is the gauge fixed action for a
theory with vanishing Hamiltonian: the first term is the gauge-fixing term,
while the second term can be written as 
\begin{equation}
  \frac{1}{2}\int_Idt\int_N \frac{1}{k!}
            \underbrace
            {\frac{(-1)^k}{(d-k)!}\vep^{i_1\ldots i_ki_{k+1}\ldots i_d}
              \dot{D}_{i_{k+1}\ldots i_d}(t,\xx)}
            _{\dot{\PHI}(t,\xx)}
            \underbrace{D_{i_1\ldots i_k}(t,\xx)}_{\PI(t,\xx)}
   \quad,
\end{equation}
which is exactly the desired expression (considering $\PHI$\ and $\PI$\ as
conjugate variables), as can be inferred from
(\ref{eq:brcc}): 
\begin{equation}
  \{
  \underbrace{D_{j_1\ldots j_k}(\xx)}_{\PI(\xx)}
  ;
  \underbrace{\frac{(-1)^k}{(d-k)!}\;\vep^{i_1\ldots i_ki_{k+1}\ldots i_d}
         D_{i_{k+1}\ldots i_d}(\yy)}_{\PHI(\yy)}
  \}=
   \de^{(d)}(\xx-\yy)\de_{i_1}^{j_1}\ldots\de_{i_k}^{j_k}
  \quad.
\end{equation}

\end{appendix}

 \vskip3em
\small
 \newcommand\wb{\,\linebreak[0]} \def\wB {$\,$\wb}
 \newcommand\Bi[1]    {\bibitem{#1}}
 \newcommand\J[5]   {{\sl #5}, {#1} {#2} ({#3}) {#4} }
 \newcommand\PhD[2]   {{\sl #2}, Ph.D.\ thesis (#1)}
 \newcommand\Prep[2]  {{\sl #2}, preprint {#1}}
 \newcommand\BOOK[4]  {{\em #1\/} ({#2}, {#3} {#4})}
 \def\jf    {J.\ Fuchs}
 \newcommand\inBO[7]  {in:\ {\em #1}, {#2}\ ({#3}, {#4} {#5}), p.\ {#6}}
 \newcommand\inBOx[7] {in:\ {\em #1}, {#2}\ ({#3}, {#4} {#5})}
 \newcommand\gxxI[2] {\inBO{GROUP21 Physical Applications and Mathematical
              Aspects of Geometry, Groups, and \A s{\rm, Vol.\,2}}
              {H.-D.\ Doebner, W.\ Scherer, and C.\ Schulte, eds.}
              \WS\Si{1997} {{#1}}{{#2}}}
 \def\anop  {Ann.\wb Phys.}
 \def\asm   {Adv.\wb Sov. \wb Math.}
 \def\comp  {Com\-mun.\wb Math.\wb Phys.}
 \def\ijma  {Int.\wb J.\wb Mod.\wb Phys.\ A}
 \def\ijmb  {Int.\wb J.\wb Mod.\wb Phys.\ B}
 \def\inma  {Invent.\wb Math.}
 \def\lmap  {Lett.\wb  Math. \wb Phys.}
 \def\mpcp  {Math.\wb Proc.\wb Camb.\wb Phil.\wb Soc.}
 \def\nupb  {Nucl.\wb Phys.\ B}
 \def\phlb  {Phys.\wb Lett.\ B}
 \def\phma  {Phil.\wb Mag.}
 \def\phrb  {Phys.\wb Rev.\ B}
 \def\phrd  {Phys.\wb Rev.\ D}
 \def\phrl  {Phys.\wb Rev.\wb Lett.}
 \def\remp  {Rev.\wb Mod.\wb Phys.} 
 \def\topo  {Topology}

 \def\AMS    {{American Mathematical Society}}
 \def\AP     {{Academic Press}}
 \def\AW     {{Addison-Wesley}}
 \def\CUP    {{Cambridge University Press}}
 \def\PUP    {{Princeton University Press}}
 \def\LMS    {{London Mathematical Society}}
 \def\MD     {{Marcel Dekker}}
 \def\NH     {{North Holland Publishing Company}}
 \def\SV     {{Sprin\-ger Ver\-lag}}
 \def\TE     {{Teubner}}
 \def\WS     {{World Scientific}}
 \def\Ad     {{Amsterdam}}
 \def\Be     {{Berlin}}
 \def\Ca     {{Cambridge}}
 \def\Le     {{Leipzig}}
 \def\NY     {{New York}}
 \def\PR     {{Providence}}
 \def\PRI    {{Princeton, New Jersey}}
 \def\RC     {{Redwood City}}
 \def\Si     {{Singapore}}


\end{document}